\documentclass[12pt,reqno,twoside]{amsart}
\usepackage{graphicx}
\usepackage{amsmath}
\usepackage{amssymb}

\usepackage{color}
\newcommand{\todo}[1]{{\color{blue} \sf To Do: [#1]}}
\newcommand{\ignore}[1]{}

\newcommand{\Q}{{\mathbb {Q}}}

\newcommand{\R}{{\mathbb{R}}}  
\newcommand{\Z}{{\mathbb{Z}}}
\newcommand{\C}{{\mathbb{C}}}

\newcommand{\N}{{\mathbb{N}}}

\newcommand\NN{\mathbb N}
\newcommand\RR{\mathbb R}
\newcommand\ZZ{\mathbb Z}
\newcommand\QQ{\mathbb Q}

\newcommand\cH{\mathcal{H}}

\newcommand\cO{\mathcal{O}}

\newcommand\norm[1]{\left\|#1\right\|}

\newcommand\abs[1]{\left|#1\right|}

\newcommand\set[1]{\left\{{#1}\right\}}

\newtheorem{thm}{Theorem}[section]
\newtheorem{lem}[thm]{Lemma}
\newtheorem{prop}[thm]{Proposition}
\newtheorem{cor}[thm]{Corollary}

\newtheorem{rem}[thm]{Remark}

\numberwithin{equation}{section}

\begin{document}

\title[Spectral gaps and intrinsic Diophantine approximation]{On discrepancy, intrinsic Diophantine approximation, and spectral gaps} 
\author{Alexander Gorodnik and Amos Nevo} 
\address{Alexander Gorodnik : Institute f\"ur Mathematik, Universit\"at Z\"urich.}
\email{alexander.gorodnik@math.uzh.ch}
\address{Amos Nevo : Department of Mathematics, Technion, and Department of Mathematics, University of Chicago.}
\email{anevo@tx.technion.ac.il, nevo@uchicago.edu}


\date{\today}



\begin{abstract} 

In the present paper we establish bounds for the size of the spectral gap for actions of algebraic groups on certain homogeneous spaces. Our approach is based on estimating operator norms of suitable averaging operators, and we develop techniques for establishing both upper and lower bounds for such norms. We shall show that this analytic problem is closely related to the arithmetic problem of establishing bounds on the discrepancy of distribution for rational points on algebraic group varieties. As an application, we show how to establish an effective bound for property $\tau$ of congruence subgroups of arithmetic lattices in algebraic groups which are forms of $SL(2)$, using estimates in intrinsic Diophantine approximation which follow from Heath-Brown’s analysis of rational points on 3-dimensional quadratic surfaces.

\end{abstract}

\maketitle

{\small \tableofcontents}
 
\section{Introduction}

\subsection{Introduction}
Let $G$ be a locally compact second countable (lcsc) group.
 A (strongly continuous) unitary representation $\pi$ of $G$ on a Hilbert space $\cH$  is said to have almost invariant vectors if there exists a sequence of unit vectors $v_n$ such that $\|\pi(g)v_n-v_n\|\to 0$
uniformly for $g$ varying in compact subsets of $G$.
If $\pi$ does not have almost invariant functions, $\pi$ is said to possess a spectral gap. This fundamental representation-theoretic property expresses the fact 
that the trivial representation of $G$ is ``isolated'' from the irreducible representations appearing in the spectral decomposition of the representation $\pi$.

The spectral gap property can be also formulated in explicit quantitative terms, via  a norm estimate of averaging
operators. For a bounded Borel measure $\beta$ on $G$, we define the bounded operator 
$$
\pi(\beta): \cH \to \cH\,\,,\,\, v \mapsto \int_G \pi(g)v\, d\beta(g).
$$
It is well-known that $\pi$ has the spectral gap property if and only if 
for some (equivalently, for all) absolutely continuous symmetric probability measure $\beta$ whose support generates  $G$, one has
$\|\pi(\beta)\|<1.$
This leads to the natural problem of establishing explicit estimates on the norm 
$\|\pi(\beta)\|$ that we aim to address and utilize in the present paper.

We will be especially interested in the spectral gap property for representations of semisimple algebraic groups acting on homogeneous spaces with finite invariant measures, which we now introduce.
Let ${\sf G}\subset \hbox{GL}_n$ be an algebraically connected linear algebraic group defined over $\Q$ which is $\Q$-simple.
We use the notations: 
$$
G_\infty:={\sf G}(\R)\quad\hbox{and}\quad
G_S:={\prod}_{p\in S} {\sf G}(\Q_p)
$$
for a finite set of primes $S$, and 
$$
\Gamma_S:={\sf G}(\Z[S^{-1}])\quad\hbox{and}\quad
\Gamma_{S,m}:=\{\gamma\in \Gamma_S:\, \gamma=e\mod m \}
$$
for $m$ coprime to $S$. We consider the spaces
$$
X_{S,m}:= (G_\infty\times G_S)/\Gamma_{S,m}
$$
and the corresponding unitary representations $\rho_{S,m}$ of the group $G_S$
acting on the spaces $L_0^2(X_{S,m})$ consisting of square-integrable functions
with zero integral.

To simplify notation,
when the set $S$ consists of a single prime $p$, 
we write $G_p$ (instead of $G_{\{p\}}$) and similarly for the other objects involving $S$ in the subscript.

 We have a  (well-defined and unique) direct integral decomposition
$$
\rho_{S,m}=\int_{\widehat G_S} \pi^{\oplus k(\pi)}\, d\Pi_{S,m}(\pi),
$$
where $\Pi_{S,m}$ is a Borel measure on the (standard Borel) unitary dual $\widehat G_S$ (consisting of equivalence classes of strongly continuous irreducible unitary representations of $G_S$) 
and $1\le k(\pi)\le \infty$ denotes multiplicities.
The support of the measure $\Pi_{S,m}$ (w.r.t. the Fell topology on the unitary dual) is called the support of the representation 
$\rho_{S,m}$.  Following \cite{BLS92}, introduce the notion of the automorphic dual of $G_S$:
$$
\widehat G_{S}^{{\rm \small aut},0}:=\overline{{\bigcup}_{(m,S)=1}\; \hbox{supp}(\rho_{S,m})}\quad\hbox{and}\quad 
\widehat G_{S}^{{\rm \small aut}} =\widehat G_{S}^{{\rm \small aut},0}\cup \{1_{G_S}\},
$$
where the closure is taken with respect to the Fell topology on the unitary dual 
$\widehat G_{S}$,
and $1_{G_S}$ denotes the trivial representation of $G_S$.
A crucially important property (called property $(\tau)$ for congruence subgroups or the Ramanjuan--Selberg property)
of the representations $\rho_{S,m}$ is that they are uniformly isolated from the trivial representation, namely, the sum of the representations $\rho_{S,m}$
has the spectral gap property. This fundamental result was proved in full generality
by Clozel \cite{Cl03}, but it has a rich history going back in some cases to works of Ramanujan and Selberg \cite{S56,S65}, Kazhdan \cite{kaz}, Burger--Sarnak \cite{BS91}, and others. We refer, for instance, to \cite{S05},\cite{C07},\cite{BB} for detailed surveys.

We introduce the height function on the group $G_S$: 
$$
\hbox{H}_S(g):={\prod}_{p\in S} \max(1,\|g_p\|_p)\quad\hbox{for $g=(g_p)_{p\in S}\in G_S$,}
$$
where $\|\cdot\|_p$ denote the maximum $p$-adic norms,
and the corresponding compact subsets
\begin{equation}\label{eq:BBh}
B_h^S:=\{g\in G_S:\, \hbox{H}_S(g)\le h \}.
\end{equation}
We consider the Haar-uniform probability measures $\beta^S_h$ supported on the sets $B^S_h$.
Then property $(\tau)$ amounts to the estimate (when $h\ge h_0(S)$)
$$
\sup \big\{\|\pi(\beta^S_h)\|:\, \pi\in\widehat G_{S}^{{\rm \small aut},0} \big\}< 1.
$$
In a number of cases, more precise explicit bounds are known.
For instance, the Ramanujan--Petersson--Selberg Conjecture amounts to showing that 
when $\sf G$ is a form of $\hbox{SL}_2$, then
$$
\sup \big\{\|\pi(\beta^S_h)\|:\, \pi\in\widehat G_{S}^{{\rm \small aut},0} \big\}\ll_{\epsilon, S} m_{G_S}(B^S_h)^{-1/2+\epsilon}\quad\hbox{for all $\epsilon>0$},
$$
where $m_{G_S}$ denotes a Haar measure on $G_S$.
It was established by Deligne in the case when $\sf G$ is anisotropic over $\R$
(see, for instance, \cite{r} for a self-contained introduction to these results), but is still open otherwise.
On the other hand, one expects only weaker decay bounds for more general algebraic groups. 

\medskip

The goal of the present paper is to develop upper and lower bounds on the norms $\|\pi(\beta^S_h)\|$ for authomorphic representations, and demonstrate their close connection and mutual relationship to problems of intrinsic Diophantine approximation on the corresponding algebraic groups. 
We note that the norm estimates are interesting only when the group $G_S$ is non-compact, which we assume throughout the paper.

We now turn to state our main results. 

\subsection{Lower bounds for operator norms in the automorphic representation}
Our first result provides an explicit lower bound for the norms of the averaging operators introduced above. 

\begin{thm}
	\label{th:lower_norm_bound_0}
	Let $\sf G$ be an algebraically connected $\Q$-simple linear algebraic group defined over $\Q$.  If $\sf L$ is a proper reductive algebraic subgroup of $\sf G$ defined over $\Q$, 
	then there exists $C > 0$ such that for all $h \ge 1$ 
	$$
	\sup\big\{ \|\pi(\beta^S_h)\|:\, \pi \in \widehat G_{S}^{{\rm\tiny aut},0} \big\}\ge  C\,  \frac{m_{L_S}(L_S\cap B^S_h)}{m_{G_S}(B^S_h)},
	$$
	where $m_{L_S}$ denotes a Haar measure on $L_S$, and the constant $C$ depends on ${\sf G},S, {\sf L}$ and the choice of Haar measures.
Therefore, for all $h \ge 1$,
$$
\sup\big\{ \|\rho_{S,m}(\beta^S_h)\|:\, m\in\N \hbox{ with } (m,S)=1 \big\}	\ge C\, \frac{m_{L_S}(L_S\cap B^S_h)}{m_{G_S}(B^S_h)}.
$$	
\end{thm}

We remark that a general, but \emph{weaker}, lower bound valid in the present context is: for all $ h \ge 1$
\begin{equation}
\label{eq:weak}
\sup\big\{ \|\pi(\beta^S_h)\|:\, \pi \in \widehat G_{S}^{{\rm\tiny aut},0} \big\}\ge  C\left(m_{G_S}(B^S_h)\right)^{-1/2}\,.
\end{equation}
We refer to Remark \ref{L-amenable} below for a justification of \eqref{eq:weak}. 
Usually, in Theorem \ref{th:lower_norm_bound_0} one can choose the subgroup 
$\sf L$, so that $m_{L_S}(L_S\cap B^S_h)\gg m_{G_S}(B^S_h)^\rho$ with $\rho>1/2$.
Then Theorem \ref{th:lower_norm_bound_0} provides a better estimate.
We exemplify this by the following explicit estimates for classical groups:

\begin{cor}\label{cor:lower}
\begin{enumerate}
	\item[1.] For ${\sf G}=\hbox{\rm SL}_n$ with $n\ge 3$,  for every prime $p$, if $h$ is sufficiently large, for all $\epsilon>0$.:
	$$
	\sup\left\{ \|\pi(\beta^{\set{p}}_h)\|:\, \pi \in \widehat G_{p}^{{\rm aut},0} \right\}\gg_\epsilon m_{G_{p}}(B^{\set{p}}_h)^{-2/n+\epsilon}\,.
	$$

	\item[2.] Let ${\sf G}=\hbox{\rm SO}_{n}$ with $n\ge 4$ and $p$ is prime 
	satisfying $p=1\,(\hbox{\rm mod}\,\,4)$. Then when $n$ is even,
	if $h$ is sufficiently large, for all $\epsilon>0$
	$$
	\sup\left\{ \|\pi(\beta^{\set{p}}_h)\|:\, \pi \in \widehat G_{p}^{{\rm aut},0} \right\}
	\gg_\epsilon m_{G_{p}}(B^{\set{p}}_h)^{-2/n+\epsilon},\;\;\; \epsilon>0,
	$$
	 and when $n$ is odd, if $h$ is sufficiently large, for all $\epsilon>0$,
	$$
	\sup\left\{ \|\pi(\beta^{\set{p}}_h)\|:\, \pi \in \widehat G_{p}^{{\rm aut},0} \right\}
	\gg_\epsilon m_{G_{p}}(B^{\set{p}}_h)^{-2/(n-1)+\epsilon}\,.
	$$

	\item[3.] Let ${\sf G}=\hbox{\rm Sp}_{2n}$ for $n\ge 2$. If $h$ is sufficiently large, for all $\epsilon>0$,
	$$
	\sup\left\{ \|\pi(\beta^{\set{p}}_h)\|:\, \pi \in \widehat G_{p}^{{\rm aut},0} \right\}
	\gg_\epsilon m_{G_{p}}(B^{\set{p}}_h)^{-2/(n+1)+\epsilon}\,.
	$$
\end{enumerate}
\end{cor}

It is to be expected that the lower bounds in Theorem \ref{th:lower_norm_bound_0} (and in some cases even the best-possible lower bounds) 
can be deduced from the description of the continuous automorphic spectrum
obtained by Langlands \cite{L76}. However our approach is different, and uses only relatively elementary considerations.

\vspace{0.2cm}

The above result raises the question of estimating the operator norms for the discrete part of the spectrum.  Let $\rho^{disc}_{S,m}$ denotes the subrepresentation of $\rho_{S,m}$
consisting of discretely embedded irreducible subrepresentations of $\rho_{S,m}$, and 
let $\rho^{cusp}_{S,m}$ denote its cuspidal subrepresentation, so that  
$\big\|\rho^{disc}_{S,m}(\beta)\big\|\ge \big\|\rho^{cusp}_{S,m}(\beta)\big\|$.    
We note that a lower bound for $\big\|\rho^{cusp}_{S,m}(\beta)\big\|$ follows from the 
remarkable equidistribution results 
for cuspidal spectrum of automorphic representations
that was recently established by 
Matz and Templier \cite{mt} and Finis and Matz \cite{fm}.   

\begin{thm}\label{th:discrete}
Let ${\sf G}$ be a split simple simply connected classical algebraic group defined over $\Q$. Let us assume that $\sf G$ is isotropic and unramified over $\Q_p$
and denote by $U_{p}$ a hyper-special maximal compact subgroup of $G_{p}$.
Then for any compactly supported $U_{p}$-biinvariant function $\beta$ on $G_{p}$,
$$
\big\|\rho^{cusp}_{p,1}(\beta)\big\|\ge \|\beta\|_{L^2}.
$$
\end{thm}

For example, when $\beta$ is a probability measure given at the normalized characteristic function of a $U_{p}$-biinvariant subset $B$ of $G_{p}$, we obtain the lower bound
$$
\big\|\rho^{disc}_{p,1}(\beta)\big\|\ge m_{G_{p}}(B)^{-1/2}\,.
$$
The bound in Theorem \ref{th:discrete} implies the bound stated in (\ref{eq:weak})
but is established for a more restricted setting, and is far less elementary compared to the discussion in Remark \ref{L-amenable}.

\subsection{Spectral gaps and mean-square discrepancy bounds on the group variety} 
Now we aim to obtain upper bounds on the averaging operators, or equivalently on the size of the spectral gap.
We shall show that such bounds can be deduced from the solution of the    
arithmetic counting problem of estimating the discrepancy of distribution of rational points.

Let now $\sf G$ be an algebraically connected simply connected $\Q$-simple linear algebraic group defined over $\Q$.
We use notations as above and additionally assume that $G_\infty$ is non-compact. 
Then the subgroup $\Gamma_{S,m}$ embedded in $G_\infty$ is dense
and, in fact, equidistributed in $G_\infty$ in a suitable sense (see \cite{GGN20,GN21}).
Let $m_{G_\infty}$ and $m_{G_S}$ denote Haar measures on $G_\infty$ and $G_S$
which are normalized so that $\Gamma_{S}$ has covolume one in 
$G_\infty\times G_S$.
We fix a left-invariant\footnote{With obvious modifications our argument applies also to a right invariant metric.}  Riemannian metric $\rho$ on $G_\infty$.
Let $B(x,r)$ denote the corresponding balls in $G_\infty$.
Let also $B$ be a bounded measurable subset of $G_S$ of positive measure.
We consider the problem of counting the points of $\Gamma_{S,m}$
contained in the domains $B(x,r)\times B$.
To study the distribution of rational points $\Gamma_{S,m}$ in $G_\infty$, 
it is natural to consider the following discrepancy function:
$$
\mathcal{D}\big(\Gamma_{S,m}, B(x,r)\times B\big):=\left| \frac{|\Gamma_{S,m} \cap (B(x,r)\times B)|}{m_{G_S}(B)} - \frac{m_{G_\infty}(B(e,r))}{|\Gamma_S/\Gamma_{S,m}|} \right|
$$
where the ball $B(x,r)$ is fixed, and the subsets $B$ eventually exhaust $G_S$.
For instance, one can take $B$ to be the height balls $B_h^S$ with $h\to\infty$.

We consider the discrepancy as a function of $x$ as it varies over $G_\infty$.
When the subset $B$ of $G_S$ is left-invariant under a compact open subgroup $U_S$ of $G_S$, the discrepancy $\mathcal{D}\big(\Gamma_{S,m}, B(x,r)\times B\big)$
is left-invariant under the subgroup $\Gamma_{S,m}\cap U_S$, which is a finite index subgroup of $\Gamma_{S,m}$.
In this case, we define the mean-square  discrepancy as 
$$
E_{S,m}(r,B):=\big\|\mathcal{D}\big(\Gamma_{S,m}, B(\cdot ,r)\times B\big)\big\|_{L^2\left((\Gamma_{S,m}\cap U_S)\backslash G_\infty)\right)},
$$
where $(\Gamma_{S,m}\cap U_S)\backslash G_\infty$ is equipped with the invariant probability measure.

Note that computing the discrepancy is a purely arithmetic problem which involves
estimating the number of rational points satisfying given Diophantine equations, inequalities and congruence conditions.
But we will show that the behavior of the discrepancy 
captures the size of the spectral gap for the corresponding automorphic representations, a purely analytic problem. The fact that the mean-square discrepancy can be used to bound the norms of the averaging operators is formulated as follows.

\begin{thm}\label{th:converse0}
	Let us assume that the set of primes $S$ is unramified
	and denote by $U_S$ a hyper-special maximal compact subgroup of $G_S$ (cf. Section \ref{sec:basic}).
	Let $B$ be any non-empty bounded $U_S$-bi-invariant subset of $G_S$ of positive measure.
	Fix an irreducible unitary representation $\pi$ of $G_S$ which is discretely embedded in $\rho_{S,m}$. Then there exists $r_0(\pi)> 0$ such the averaging operator $\pi(\beta)$ supported on $B$ satisfies the bound
	\begin{equation}\label{upper-bound}
	 \norm{\pi(\beta)}\ll\, \abs{\Gamma_S/\Gamma_{S,m}} r^{-\dim(G_\infty)} E_{S,m}(r,B)
	\end{equation}
	for all $0< r\le r_0(\pi)$. Here the 
	implied constant depends only on $G_\infty$.
\end{thm}

While our approach has little in common with the work of Sarnak and Xue \cite{sx}, 
their work was the main inspiration for Theorem \ref{th:converse0}, which  
also develops an approach for establishing the uniform spectral gap property 
for the Archemedian factor based on lattice point counting estimates.

Let us note the following regarding Theorem \ref{th:converse0}:
\begin{itemize}
\item For an upper bound on $\norm{\pi(\beta)}$ to be meaningful it must be strictly less that $1$. As we shall see below, in suitable circumstances this is indeed the case for the bound (\ref{upper-bound}), 
provided that the measure of the set $B$ is sufficiently large compared with the inverse of the measure of the ball $B(e,r)$. 
\item The assumption that $\pi$ is discretely embedded in 
$\rho_{S,m}$ is important for our present argument.
We expect however that this argument can be developed further to 
deal with all representations which are weakly contained in 
 $\rho_{S,m}$. Note that the complete description of continuous spectrum
is known due to the work of Langlands \cite{L76}. Therefore, understanding the discretely embedded representations is 
the crucial missing ingredient. 
\end{itemize}

In our previous work \cite{GGN20,GN21}, we have shown that analysis on the homogeneous spaces
$X_{S,m}$ can be used to establish mean, almost everywhere, and pointwise bounds
for the discrepancy function. Remarkably, Theorem \ref{th:converse0} (and more generally Theorem \ref{th:converse} below) show that 
the converse is also true, and estimates on discrepancy lead to quantitative 
estimates for the spectral gap. 
These observations raise an interesting possibility of studying 
the spectral gap problem via arithmetic consideration.
Here we carry this out for forms of $\hbox{SL}_2$.

\subsection{Property $(\tau)$ for congruence subgroups}
Let $\sf G$ be a linear algebraic group defined over $\Q$ which is a form 
of $\hbox{SL}_2$. More explicitly, $\sf G$ can be viewed as the set
of norm one elements in a division algebra defined over $\Q$. 
The integral structure on $\sf G$ is defined with respect to an order of the division algebra. We fix a prime $p$, and assume that $\sf G$ is isotropic over $\R$ and over $\Q_p$. Moreover, we assume that ${\sf G}(\Z_p)$ is a hyper-special maximal compact subgroup of $G_{p}$ (which is the case for all but finitely many primes).
In this setting, we derive a bound on the norms of automorphic representations: 

\begin{thm}\label{th:sl20}
For every $\ell\in \mathbb{N}$ coprime to a prime $p$ satisfying the foregoing conditions,
	$$
	\big\|\rho_{p,\ell}(\beta_h^{\set{p}})\big\|\ll_{p,\ell} m_{G_{p}}\big(B_h^{\set{p}}\big)^{-\kappa},
	$$
	where $\kappa=1/4$ if ${\sf G}$ is anisotropic over $\Q$,
	and $\kappa=1/16$ if ${\sf G}$ is isotropic over $\Q$.
	As a result, for any $\pi\in \widehat{G}_{p}^{{\rm\tiny aut},0}$,
	$$
     \big\|\pi(\beta_h^{\set{p}})\big\|\ll_{p,\ell} m_{G_{p}}\big(B_h^{\set{p}}\big)^{-\kappa}.
    $$
\end{thm}

Our estimate falls short of the best known bound, which corresponds to $\kappa={25}/{64}$. This bound is due to Kim and Sarnak \cite[Appendix 2]{ks} over $\Q$ (and Blomer and Brumley \cite{bb11} over number fields) and was proved by quite a different argument. Our approach involves establishing a bound on the discrepancy
function directly using the refined circle method arguments due to Heath-Brown, 
and then combining this bound with a version of Theorem \ref{th:converse0}.
We remark that for technical reasons, we will work with smooth test-functions on $G_\infty$ rather than the balls as above.

\subsection{Organisation of the paper}
In the next section we set up notation and review basic facts that will be used throughout the paper.
Then in Section \ref{sec:lower} we develop a method for proving 
lower bounds for the norms of averaging operators and establish 
Theorem \ref{th:lower_norm_bound_0} and Corollary \ref{cor:lower}.
In Section \ref{sec:descrete} we estimate the norms of averaging operators
acting on the discrete part of the spectrum and prove Theorem \ref{th:discrete}.
In Section \ref{sec:spectral bounds} we show how to derive norm bounds from estimates
on the discrepancy of rational points and prove Theorem \ref{th:converse0}.
This approach will be utilized in Section \ref{sec:sl2}
where we use the refined circle method technique developed by Heath-Brown \cite{HB96}  to establish the required estimates 
on discrepancy and prove Theorem \ref{th:sl20}.

\subsection*{Acknowledgements}
A.G. was supported by SNF grant 200021--182089. 
A.N. was partially supported by ISF Moked Grant 2919-19.

\section{Basic notation}\label{sec:basic}
	
	Let ${\sf G}\subset {\sf G}{\sf L}_n$ be a connected semisimple linear algebraic group defined over a number field $K$. We denote by $K_v$, $v\in V$ the completions of $K$.
	We fix a finite set $S$ of non-Archemedian completions such that 
	the group $\sf G$ is isotropic for at least one completion from $S$.
	We denote by $G_\infty$ the product of ${\sf G}(K_v)$ over the Archemedian 
	completions and by $G_S$  the product of ${\sf G}(K_v)$ over $S$.
	We consider $G_\infty$ and $G_S$ as locally compact groups 
	equipped with the topology arising from from the field completions.
	We fix a choice of Haar measures $m_{G_\infty}$ and $m_{G_S}$ on $G_\infty$ and $G_S$ respectively.

	Let $O_S:=\{x\in K:\, |x|_v\le 1\hbox{ for $v\notin S$}\}$
	denote the subring of $K$ consisting of elements which are $K_v$-integral for $v$ outside $S$ (also known as the ring of $S$-integers in $K$). 
	For non-Archemedian completions, we write $O_v$ for the ring of integers in $K_v$ (note that this notation means that $O_v\neq O_{\set{v}}$ !). 
	We denote by $\Gamma_S:={\sf G}(O_S)$ the $S$-arithmetic subgroup of ${\sf G}(K)$.
	We denote by $\mathcal{O}^S$ the product of ${\sf G}(O_v)$ over 
	the non-Archemedian completions $v$ which are not in $S$.
	Then $\Gamma_S$ can be viewed as a subgroup $\mathcal{O}^S$,
	and for an open subgroup $\mathcal{O}$ of $\mathcal{O}^S$,
	we define $\Gamma_{S,\mathcal{O}}:=\Gamma_S\cap \mathcal{O}$.
	Under the natural embedding, we view $\Gamma_{S,\mathcal{O}}$
	as a subgroup of $G_\infty\times G_S$. Then it is a discrete subgroup with finite covolume, and we consider the homogeneous spaces 
	$X_{S,\mathcal{O}}:=(G_\infty\times G_S)/\Gamma_{S,\mathcal{O}}$
	equipped with the invariant probability measure $\mu_{S,\cO}$.
	Let $\rho_{S,\mathcal{O}}$ denote the corresponding unitary representations
	of the group $G_S$ on the spaces
	$L_0^2(X_{S,\mathcal{O}})$, the space of $L^2$-integrable functions with zero integral. 
	
	The group $\sf G$ is said to be unramified over a completion $K_v$ 
	if it is quasi-split and split over an unramified extension of $K_v$.
	If $\sf G$ is unramified over $K_v$, there exists a canonical conjugacy class 
	of maximal compact open subgroups of $G_v$ --- the so-called hyperspecial 
	maximal compact subgroups (see \cite{tits}). 
	We fix a choice $U_v$ of a hyperspecial 
	maximal compact subgroup $U_v$ of $G_v$. More generally, we write $U_S:=\prod_{v\in S} U_v$ when $\sf G$ is unramified over every $v \in S$ (i.e. unramified over $S$). The pair $(G_S,U_S)$ is known to have the Gelfand-Selberg property, namely the convolution algebra $L^1(U_S\backslash G_S/U_S)$
	is commutative. We say that a unitary representation of $G_S$ is spherical
	if there exists a non-zero $U_S$-invariant vector. 
	It follows from the Gelfand-Selberg property that if a spherical representation is irreducible, then 
	the dimension of the space of $U_S$-invariant vectors is at most one. 
		
	Let $\widehat G_S$ denote the unitary dual of the group $G_S$
	equipped with the Fell topology. We recall that given an arbitrary unitary representation 
	$\sigma$ of  $G_S$ on a Hilbert space $\mathcal{H}_\sigma$, we have the direct integral decomposition
	$$
	\sigma=\int_{\widehat G_S} \pi^{\oplus k(\pi)}\, d\Pi_{S,\mathcal{O}}(\pi)
	$$
	where $\Pi_{S,\mathcal{O}}$ is a Borel measure on $\widehat G_S$, and $1\le k(\pi)\le \infty$ denotes 
	the multiplicities. 
	The support  of the measure $\Pi_{S,\mathcal{O}}$ is uniquely defined and is called 
	the support  of the representation $\sigma$, denoted $\hbox{supp}(\sigma)$.
	One says that $\sigma$ is weakly contained in a subset $\Omega$ of $\widehat G_S$ 
	if $\hbox{supp}(\sigma)\subset \overline{\Omega}$,
	where the closure is taken with respect to the Fell topology. 
	Namely, this means that every function of positive type associated to $\sigma$
	can be approximated, uniformly on compact subsets of $G_S$, by finite sums of functions
	of positive type associated to $\pi\in\Omega$. 
	More generally, for unitary representations $\sigma_1$ and $\sigma_2$ of $G_S$, we say that $\sigma_1$ is weakly contained in $\sigma_2$ 
	if $\hbox{supp}(\sigma_1)\subset \hbox{supp}(\sigma_2)$.
	We refer to \cite[App.~F]{Kazhdan} for basic properties of the notion of weak containment.
	In particular, we recall that $\sigma_1$ is weakly contained in $\sigma_2$ is equivalent to
	\begin{equation}
	\label{eq:norm0}
	\|\sigma_1(F)\|\le \|\sigma_2(F)\|\quad\hbox{ for all $F\in L^1(G_S)$,}
	\end{equation}
	and that $\sigma$ is weakly equivalent to the representation $\oplus\{\pi:\, \pi\in\hbox{supp}(\sigma)\}$. We say that a representation is tempered 
	if it is weakly contained in the regular representation $L^2(G_S)$.
	
	Following \cite{BLS92}, we introduce the notion of automorphic dual of $G$ over $S$:
	$$
	\widehat G_{S}^{{\rm \small aut}}:=\widehat G_{S}^{{\rm \small aut},0}\cup \{1_{G_S}\}
	\quad\hbox{with}\quad \widehat G_{S}^{{\rm \small aut},0}:= \overline{{\bigcup}_{\mathcal{O}\subset \mathcal{O}^S}\; \hbox{supp}(\rho_{S,\cO})},
	$$
	where the union is taken over open subgroups
	$\mathcal{O}$ of $\mathcal{O}^S$, and 
	the closure is with respect to the Fell topology on $\widehat G_S$.
Given a subset $B$ of $G_S$ with finite positive measure, 
we denote by $\beta$ the Haar-uniform probability measure supported on the set $B$.
For $\pi\in\widehat G_S$, we consider the averaging operators $\pi(\beta)$ which are explicitly defined by
$$
\pi(\beta)=\frac{1}{m_S(B)}\int_{B} \pi(g)\, dm_S(g).
$$
Recall that for any two finite Borel measures $\beta_1$ and $\beta_2$ on $G_S$
$$
\pi(\beta_1*\beta_2)=\pi(\beta_1)\pi(\beta_2)\quad\hbox{and}\quad \pi(\beta)^*=\pi(\beta^*).
$$
We will be especially interested in the sets 
$$
B_h^S:=\{g\in G_S:\, \hbox{H}_S(g)\le h \}
$$
defined by the height function
$$
\hbox{H}_S(g):=\prod_{v\in S} \max(1,\|g_v\|_v)\quad\hbox{for $g=(g_v)_{v\in S}\in G_S$.}
$$
In this case, we denote the corresponding probability measure by $\beta_h^S$ (and when no ambiguity is present, by $\beta_h$).
	
For some classes of unitary representations, explicit estimate on the operator 
norms $\|\pi(\beta)\|$ have been established. For instance, according to the Kunze--Stein
Phenomenon (proved for $p$-adic groups in \cite{ve} following the method of \cite{c78}), for every $q\in [1,2)$,
$$
\|f*\phi\|_2\ll_q \|f\|_q \|\phi\|_2\quad\quad\hbox{for all $f\in L^q(G_S)$ and $\phi\in L^2(G_S)$.}
$$
It follows (see \cite{GN15} for a full account)  that for every tempered representation $\pi$, and for every subset $B\subset G_S$ of positive finite measure 
\begin{equation}
\label{eq:temp}
\|\pi(\beta)\|\ll_{\epsilon} m_S(B)^{-1/2+\epsilon}
\end{equation}
for all $\epsilon>0$. More generally, if the tensor power $\pi^{\otimes n}$
is tempered for some even integer $n$, then by \cite{N2}, \cite{GN10} the following bound holds :
$$
\|\pi(\beta)\|\ll_{\epsilon} m_S(B)^{-1/(2n)+\epsilon}
$$
for all $\epsilon>0$.

In particular, the tensor power property is known to hold for non-trivial irreducible unitary representations of $\hbox{SL}_d(K_v)$ (and more generally simple groups with property $T$). Therefore, for every such representation $\pi$,
$$
\|\pi(\beta)\|\ll m_S(B)^{-\kappa(\pi)}
$$
for some exponent $\kappa(\pi)>0$, which is uniform in the representation $\pi$ when $d>2$. 

When an irreducible representation $\pi$ is spherical, then 
\begin{equation}
\label{eq:norm_spher}
\norm{\pi(\beta)}^2 =\norm{\pi(\beta^\ast \ast \beta)}=\abs{\frac{1}{m_S(B)}\int_{B} \omega_\pi(g)\, dm_S(g)}^2\le \frac{1}{m_S(B)}\int_{B} \abs{\omega_\pi(g)}^2\, dm_S(g),
\end{equation}
where $\omega_\pi$ denotes the spherical functions associated to the representation $\pi$. 
Then $\|\pi(\beta)\|$ can be analyzed using estimates on the spherical functions. For instance, considering the complementary series representations
of $\hbox{SL}_2(\Q_p)$, one can observe that the norms $\|\pi(\beta)\|$
may decay arbitrary slowly for irreducible representations $\pi$ and thus also construct 
examples of (reducible) representations without invariant vectors such that $\|\pi(\beta)\|=1$.

\section{Lower bounds for norms of automorphic representations}\label{sec:lower}

Let $\sf G$ be a connected semisimple linear algebraic group defined over a number field $K$. We keep the notation and assumptions introduced in Section \ref{sec:basic}.
Our goal is to establish a lower bound for $\|\pi(\beta^S_h)\|$, for the representations
$\pi\in \widehat G_{S}^{{\rm \small aut},0}$.

\begin{thm}\label{lower-bound}
	\label{th:lower_norm_bound}
	Let $\sf L$ be \emph{any} proper reductive algebraic subgroup of $\sf G$ defined over $K$.  If $L_S$ is non-amenable, then for some positive constant $C$
	 and for all $h\ge 1$,  
	$$
	\sup\left\{ \|\pi(\beta^S_h)\|:\, \pi \in \widehat G_{S}^{{\rm\tiny aut},0} \right\}	\ge C  \frac{m_{L_S}(L_S\cap B^S_h)}{m_{G_S}(B^S_h)},
	$$
where $C$ depends on ${\sf G},S, {\sf L}$ and the choice of Haar measures, but not on $h$.

 In particular,
for all $h\ge 1$,
$$
\sup\big\{ \|\rho_{S,\mathcal{O}}(\beta^S_h)\|:\, \hbox{$\mathcal{O}$ --
	an open subgroup of $\mathcal{O}^S$}\big\} \ge C  \frac{m_{L_S}(L_S\cap B^S_h)}{m_{G_S}(B^S_h)}.
$$
\end{thm}
 
\begin{proof}
The key starting point of our argument are the fundamental results from \cite{BLS92}
and \cite[Sec.~3.3]{C07} which assert that if $\pi\in \widehat L_S^{\rm\tiny aut}$, then 
the induced representation $\hbox{Ind}_{L_S}^{G_S}(\pi)$ is weakly contained in
$\widehat G_S^{\rm\tiny aut}$.
In particular, if we denote by $\tau$ the unitary representation of $G_S$ on $L^2(G_S/L_S)$ (namely the representation induced from the trivial representation of $L_S$),
then it follows that $\tau$ is weakly contained in $\widehat G_S^{\rm\tiny aut}$.

We also note that the representation $\tau$ does not weakly contain 
the trivial representation $1_{G_S}$. Indeed, consider any 
lcsc group $G$ and a closed subgroup $L$, with $G$ acting on the homogeneous space $G/L$ with infinite invariant measure. Then weak containment of the trivial representation of $G$ in $L^2(G/L)$ is equivalent to the existence of an asymptotically invariant (=F\o lner) sequence of subsets of $G/L$. The proof this statement is due to Eymard  \cite{Ey72}  and Greenleaf \cite[Thm. 4.1]{Gr71}. We refer to \cite[Thm. 3.3]{IN96} for a short proof, and to \cite[\S \S 3.1,4.1]{AD03} for a full account of this topic.  

Furthermore, assume that $G$ is the $K$-points of an algebraically connected semisimple algebraic group over a locally compact non-discrete field $K$ of characteristic zero. Then  the existence of an asymptotically invariant sequence on $G/L$ implies that $L$ is Zariski-dense in $G$. This follows immediately from \cite[Thm. 6.1]{IN96}, and thus when $L$ is a proper unimodular algebraic subgroup, we conclude that the trivial representation is not weakly contained in $L^2 (G/L)$. 
This argument also extends to the case of finite products of algebraic groups over local fields of characteristic zero.

We can therefore conclude that $\tau$ is in fact weakly contained in $\widehat G_S^{{\rm\tiny aut},0}$.
In particular, 
$\tau$ is weakly contained in the direct sum of $\pi$ in $\widehat G_S^{{\rm\tiny aut},0}$ (see \cite[Prop.~F.2.7]{Kazhdan}).
Therefore, we deduce that
\begin{equation}
\label{eq:norm}
\sup\left\{ \|\pi(\beta_h^S)\|:\, \pi \in \widehat G_{S}^{{\rm\tiny aut},0} \right\}\ge \|\tau(\beta^S_h)\|.
\end{equation}
Hence, to finish the proof, it remains to establish a suitable lower bound for 
$\|\tau(\beta^S_h)\|$.

Let $\nu$ be a $G_S$-invariant measure on the homogeneous space $G_S/L_S$.
We denote by $p:G_S\to G_S/L_S$ the factor map.
Let $\phi$ be a characteristic function a compact neighborhood $U$ of the identity coset in $G_S/L_S$. Then  
\begin{align*}
\left<\tau(\beta^S_h)\phi,\phi\right>=\frac{1}{m_{G_S}(B^S_h)}\int_{G_S/L_S} \int_{B^S_h} \phi(g^{-1}xL_S)\phi(xL_S)\, dm_{G_S}(g)d\nu(xL_S).
\end{align*}
We observe that $\phi(g^{-1}xL_S)=1$ provided that $g\in xp^{-1}(U)^{-1}$,
and it is equal to zero otherwise. Hence,
\begin{align*}
\left<\tau(\beta^S_h)\phi,\phi\right>=\frac{1}{m_{G_S}(B^S_h)}\int_{U} 
m_S\big(B^S_h\cap xp^{-1}(U)^{-1}\big) \, d\nu(xL_S).
\end{align*}
If $x$ belongs to a fixed compact subset of $G_S$, then 
it follows from properties of $p$-adic norms that 
there exists $c_1>0$ such that $\hbox{H}_S(xg)\le c_1\, \hbox{H}_S(g)$ for all $g\in G_S$.
Therefore $x^{-1}B^S_h\supset B^S_{c_1^{-1} h}$, and we deduce that 
$$
m_{G_S}\big(B^S_h\cap xp^{-1}(U)^{-1}\big)=m_{G_S}\big(x^{-1} B^S_h\cap p^{-1}(U)^{-1}\big)\ge m_{G_S}\big(B^S_{c_1^{-1} h}\cap p^{-1}(U)^{-1}\big),
$$
and 
\begin{align*}
\left<\tau(\beta^S_h)\phi,\phi\right>\ge \nu(U) \frac{m_{G_S}\big(B^S_{c_1^{-1}h}\cap p^{-1}(U)^{-1}\big)}{m_{G_S}(B^S_h)}.
\end{align*}
If the neighbourhood $U$ is chosen to be sufficiently small, there exists a  continuous section $\sigma:U\to p^{-1}(U)$ of the factor map $p$ such that the map $L_S\times U \to p^{-1}(U)^{-1}$ defined by
$(u,l)\mapsto l\sigma(u)^{-1}$ is a homeomorphism.
Since $\sigma(U)$ is compact, there exists $c_2>0$ such that
$\hbox{H}_S(gx)\le c_2\, \hbox{H}_S(g)$ for all $g\in G_S$ and $x\in \sigma(U)^{-1}$,
so that
$$
B^S_{c_1^{-1} h}\cap p^{-1}(U)^{-1}=B^S_{c_1^{-1}h}\cap L_S \sigma(U)^{-1}\supset (B^S_{c_1^{-1}c_2^{-1}h}\cap L_S)\sigma(U)^{-1}.
$$
We observe that the Haar measure on $p^{-1}(U)\subset G_S$ can be decomposed as
$$
\int_{G_S} f\, dm_S=\int_{U} \int_{L_S} f(l\sigma(u)^{-1})\, dm_{L_S}(l) d\nu(u),
$$ 
so that 
$$
m_{G_S}\big((L_S\cap B^S_{c_1^{-1}c_2^{-1}h})\sigma(U)^{-1}\big)\ge \nu(U) m_{L_S}\big(L_S\cap B^S_{c_1^{-1}c_2^{-1}h}\big).
$$
Therefore, combining the above estimates, we conclude that 
\begin{align*}
\left<\tau(\beta^S_h)\phi,\phi\right>\gg \frac{m_{L_S}(L_S\cap B^S_{c_1^{-1}c_2^{-1}h})}{m_{G_S}(B^S_h)}.
\end{align*}
Hence, it follows from the definition of the operator norm that 
\begin{align*}
\|\tau(\beta^S_h)\|\gg \frac{m_{L_S}\big(L_S\cap B^S_{c_1^{-1}c_2^{-1}h}\big)}{m_{G_S}(B^S_h)}.
\end{align*}
To conclude the proof of the first estimate stated in Theorem \ref{lower-bound}, we recall \eqref{eq:norm}
and use a property of the volume function formulated and established in Lemma \ref{l:vol} immediately below.

To prove the second inequality stated in Theorem \ref{lower-bound}, we observe that every $\pi \in \widehat G_{S}^{{\rm\tiny aut},0}$ is weakly contained in
$\oplus_{\mathcal{O}}\{\rho_{S,\mathcal{O}}\}$. 
Hence, it follows from \eqref{eq:norm0} that
$$
\|\pi(\beta^S_h)\|\le {\sup}_\mathcal{O}\; \|\rho_{S,\mathcal{O}}(\beta^S_h)\|,
$$
and the second inequality follows directly from the first estimate.
\end{proof}

\begin{lem}\label{vol-estimate}
	\label{l:vol}
	Let ${\sf L}\subset \hbox{\rm GL}_n$ be a linear reductive algebraic group defined over a number field $K$
	and $S$ a finite set of non-Archemedean completions of $K$.
	We denote by $m_{L_S}$ a Haar measure on the group $L_S$.
	Then there exists $c>0$ such that the sets 
	$B_h^S\cap L_S:=\{g\in L_S:\, \hbox{\rm H}_S(g)\le h\}$ satisfy
	$$m_{L_S}(B^S_{2h}\cap L_S )\le c\, m_{L_S}(B^S_h\cap L_S)$$
	for all sufficiently large $h$.
\end{lem}

\begin{proof}
First, we consider the case when $S=\set{v}$ consist of a single completion $v$. For notational simplicity, in the present proof we set $W(K_v)=W_v$ (rather than $W_{\set{v}}$) for any algebraic group $W$. 
Let us recall the Cartan decomposition of $L_v$ (see, for instance, \cite[Ch.~0]{Si}).
We take a maximal $K_v$-split torus $\sf T$ of $\sf L$ and $\sf P$ a minimal parabolic 
subgroup associated to $\sf T$ which has a decomposition $\sf P=MU$, where $\sf M$ is the centralizer of $\sf T$ in $\sf L$, and $\sf U$ is the unipotent radical.
We denote by $\Sigma^+$ the set of positive roots of $\sf A$
associated to the parabolic subgroup $\sf P$ and $\Pi\subset \Sigma^+$
the set of simple roots.
We write $\mathcal{X}(\sf M)$ for the group of algebraic characters of $\sf M$. 
Given $\chi\in \mathcal{X}(\sf M)$, we write $|\chi(m)|_v=q_v^{\left<\chi,\omega(m) \right>}$ for $m\in M_v$, where $\omega: M_v\to \hbox{Hom}(\mathcal{X}(\sf M),\mathbb{Z})$. We denote by $M_v^0$ the kernel of $\omega$.
Then $M_v/M_v^0$ is free abelian group, and moreover $\omega(M_v)$ can be considered as a lattice in $\hbox{Hom}(\mathcal{X}(\sf M),\mathbb{R})$.
We set 
$$
M_v^+:=\{m\in M_v:\, \left<\alpha,\omega(m)\right>\ge 0\hbox{ for all $\alpha\in \Sigma^+$}  \}.
$$
Let $U_v$ be a good maximal compact subgroup of $L_v$
associated to $A_v$. Then the Cartan decomposition holds
$$
L_v={\bigsqcup}_{m\in M_v^+/M_v^0} U_v m U_v.
$$
We also set $A^+_v:=A_v\cap M_v^+$ and $A^0_v:=A_v\cap M_v^0$.
Then this decomposition can be rewritten as 
$$
L_v={\bigsqcup}_{a\in A_v^+/A_v^0,\, \omega\in\Omega_v} U_v a\omega U_v,
$$
where $\Omega_v$ is a finite subset of $M_v$.

Let us consider the representation $L_v\to \hbox{GL}_n(K_v)$.
Since the torus $\sf T$ is $K_v$-split, the image of $T_v$ is diagonalizable,
and we denote by $\Phi\subset \mathcal{X}(\sf A)$ the set of the corresponding weights. We introduce a modified height function $\hbox{H}'_v:L_v\to\mathbb{R}^+$
defined by 
$$
\hbox{H}'_v(g):=\max(1,\|a\|'_v) \quad\hbox{when $g\in U_v a\omega U_v$},
$$
where $\|a\|'_v :=\max_{\chi\in \Phi} |\chi(a)|_v$.
It follows from the basic properties of norms and compactness that
there exist $c_1,c_2>0$ such that 
$c_1\,\hbox{H}'_v \le \hbox{H}_v \le c_2\,\hbox{H}_v'$.
Therefore, it will be sufficient to analyze measures of the sets 
$B_h^\prime:=\{g\in L_v:\, \hbox{\rm H}'_v(g)\le h\}$.
We obtain
$$
m_{L_v}(B_h^\prime)={\sum}_{a\in A_v^+/A_v^0: \hbox{\rm\tiny H}'_v(a)\le h,\,\, \omega\in\Omega_v}  m_{L_v}(U_v a\omega U_v).
$$
Since $\omega U_v\omega^{-1}\cap U_v$ has finite index in both $U_v$ and $\omega U_v\omega^{-1}$, it is clear that 
$$
m_{L_v}(U_v a U_v)\ll m_{L_v}(U_v a\omega U_v)\ll m_{L_v}(U_v aU_v).
$$
Let $\delta:M_v\to\mathbb{R}^+$ denote the modular function of $P_v$.
Then according to \cite[Lemma~4.1.1]{Si}, there exist $c_1,c_2>0$ such 
that 
$$
c_1\, \delta(m)\le m_S(U_v m U_v)\le c_2\, \delta(m)\quad\hbox{for every $m\in M_v^+$.}
$$
Hence, it remains to investigate the function
$$
V(h):={\sum}_{a\in A_v^+/A_v^0: \hbox{\rm\tiny H}'_v(a)\le h}  \delta(a).
$$
Let us fix a basis of $A_v/A_v^0$ which is dual to the basis $\Pi$ for simple roots.
Then $A_v/A_v^0\otimes \R$ can be identified with $\R^r$, so that 
$\Lambda_v:=A_v/A_v^0$ is a lattice in $\R^r$.
For $\chi\in \Phi$, $|\chi(a)|_v=q_v^{\sum_{i=1}^r n_i(\chi)t_i}$,
where $n_i(\chi)\in \Q$ and $(t_1,\ldots,t_r)$ denote the coordinates of $a$.
Similarly, $\delta(a)=q_v^{\sum_{i=1}^r m_i t_i}$ with $m_i\in\N$.
Using this notation, we rewrite $V(h)$ as 
$$
V(h):={\sum}_{t\in \Lambda_v \cap D_h}  q_v^{\sum_i m_i t_i},
$$
where 
$$
D_h:=\left\{t\in (\R^+)^r:\,  \left(\sum_{i=1}^r n_i(\chi)t_i\right)^+ \le \log_{q_v} (h)\quad\hbox{ for $\chi\in \Phi$} \right\}.
$$
Here we use the notation $x^+:=\max(0,x)$.
We shall also consider the integral
\begin{equation}
\label{eq:iii}
I(h):= \int_{D_h} q_v^{\sum_i m_i t_i}\, dt.
\end{equation}
Comparing $I(h)$ with a suitable Riemann sums for $I(h)$, we deduce that that there exist $c_1,c_2>0$ and $d_1,d_2>0$ such that
$$
c_1\, I(h-d_1)\le V(h)\le c_2\, I(h+d_2)
$$
for all sufficiently large $h$.
Finally, it follows by a change of variables that there exists $c>0$ such that 
$I(2h)\le c\, I(h)$ for $h\ge 1$. Hence, the similar estimate also holds for $V(h)$, which completes the proof of the lemma when $S$ consists of a single completion.

For general $S$, we observe that the set $B^S_h\cap L$ is defined by 
the condition $\sum_{v\in S} \log \hbox{H}_v(g_v)\le \log h$.
Therefore the required estimate follows from \cite[Prop.~7.7]{GN10}.
\end{proof}

\begin{rem}\label{L-amenable}
{\rm
Here we outline a proof of the estimate \eqref{eq:weak}.
It follows from the the proof of Theorem \ref{lower-bound} (with $L_S=\{e\}$)
that the regular representation
$\lambda_{G_S}$ is weakly contained in the automorphic spectrum
$\widehat G_{S}^{{\rm\tiny aut},0}$.  Hence, for any absolutely continuous probability measure $\beta$ 
$$\norm{\lambda_{G_S}(\beta)}\le \sup\left\{ \|\pi(\beta)\|:\, \pi \in \widehat G_{S}^{{\rm\tiny aut},0} \right\}\,.$$ 
Let $P_S$ be a minimal parabolic subgroup of $G_S$.
Since $P_S$ is amenable, the regular representation $\lambda_{P_S}$ weakly contains trivial representation $1_{P_S}$. Therefore, the induced representation $\lambda_{G_S}=\hbox{Ind}_{P_S}^{G_S}(\lambda_{P_S})$ weakly contains 
the induced representation $\sigma:=\hbox{Ind}_{P_S}^{G_S}(1_{P_S})$ (see \cite[F.3.5]{Kazhdan}), and in particular,
$$
\norm{\sigma(\beta)}\le \norm{\lambda_{G_S}(\beta)}.
$$
This suggests that the Harish-Chandra function $\Xi_{G_S}$, which is the spherical function associated to the representation $\sigma$ can be used 
to estimate the norm from below. Indeed, it follows from \eqref{eq:norm_spher}
that whenever the measure $\beta$ is bi-invariant under $U_S$, 
$$ \norm{\sigma(\beta)}=\int_{G_S} \Xi_{G_S}(g)d\beta(g)$$ 
since $\Xi_{G_S}\ge 0$.
We recall furthermore that the  Harish-Chandra function satisfies the inequality 
 $$ \Xi_{G_S}(g)\ge  C^\prime\,\delta^{1/2}(a(g)),$$ 
 where $\delta(g)$ is the modular character of a minimal parabolic subgroup, and $a(g)$ is the Cartan component of $g$ w.r.t. the associated Cartan decomposition. 
 (See e.g. \cite{Si} or \cite[Prop. 2.1]{Ta83} in the totally disconnected case, and \cite[Thm. 4.6.5]{GV88} in the real case). Therefore, 
 $$ \norm{\sigma(\beta_h^S)}\ge C^\prime\,  m_{G_S}(B_h^S)^{-1}\int_{B_h^S}\delta^{1/2}(a(g))\, dm_{G_S}(g).$$
The argument of the proof of Lemma \ref{vol-estimate} applied to $G_S$  
and using the analysis of integrals of the form \eqref{eq:iii}, gives that 
 $$ \int_{B_h^S}\delta^{1/2}(a(g))\, dm_{G_S}(g)\ge C''\, m_{G_S}(B_h^S)^{1/2}\,.$$
 This implies \eqref{eq:weak}.
} 
\end{rem}

\smallskip

We now turn to Corollary \ref{cor:lower}, and note that the proof of Theorem  \ref{th:lower_norm_bound} also provides a method for estimating the volumes $m_S(B^S_h)$. We demonstrate the results by executing explicit computations for some classical groups.

\begin{proof}[Proof of Corollary \ref{cor:lower}]
We now fix $S=\set{p}$. Since we know that
\begin{equation}
\label{eq:est}
\sup\left\{ \|\pi(\beta^{\set{p}}_h)\|:\, \pi \in \widehat G_{p}^{{\rm\tiny aut},0} \right\}	\gg  \frac{m_{L_{p}}(L_{p}\cap B^{\set{p}}_h)}{m_{G_{p}}(B^{\set{p}}_h)},
\end{equation}
it remains to estimate the relevant volumes. 
It follows from the proof of Lemma \ref{l:vol} that 
$$
I(h)\ll m_{G_{p}}(B^{\set{p}}_h)\ll I(h),
$$ where $I(h)$ is the integral defined in \eqref{eq:iii}.
We denote by $\alpha$ the maximum of $\sum_{i=1}^r m_i t_i$ on the domain $D_1$.
We observe that 
$$
h^\alpha\ll I(h)\ll \hbox{vol}(D_1) h^\alpha\ll_{\epsilon} h^{\alpha+\epsilon}
$$
for all $\epsilon>0$. In all the cases considered, the corresponding adjoint representations on the Lie algebra of $\sf G$ 
are irreducible, and the domain $D_1$ is given by
$$
D_1:=\left\{t\in (\R^+)^r:\,  \sum_{i=1}^r n_i t_i \le 1 \right\},
$$
where $\sum_{i=1}^r n_i t_i$ is the highest weight of the representation.
Then $\alpha=\max(m_i/n_i)$.
These considerations can be used to estimate the norms for classical groups.

For ${\sf G}=\hbox{SL}_n$ with $n\ge 3$, we apply the estimate \eqref{eq:est}
with the subgroup ${\sf L}=\hbox{SL}_{n-1}$. We obtain:
$$
m_{G_{p}}(B^{\set{p}}_h)\gg h^{n^2-n}\quad\hbox{and} \quad m_{L_{p}}(L_{p}\cap B^{\set{p}}_h)\ll_\epsilon h^{(n-1)^2-n+1+\epsilon}
$$
for all $\epsilon>0$, and it follows that
$$
\sup\left\{ \|\pi(\beta^{\set{p}}_h)\|:\, \pi \in \widehat G_{p}^{{\rm\tiny aut},0} \right\}	\gg_\epsilon h^{-2(n-1)+\epsilon } \gg m_{G_{p}}(B^{\set{p}}_h)^{-2/n+\epsilon}.
$$

Let ${\sf G}=\hbox{SO}_{n}$ with $n\ge 4$. The assumption that $p=1\,(\hbox{mod}\, 4)$
implies that  $\sf G$ is split over $\Q_p$.
We apply the estimate \eqref{eq:est}
with the subgroup ${\sf L}=\hbox{SO}_{n-1}$.
When $n$ is even, we have
	$$
	m_{G_{p}}(B^{\set{p}}_h)\gg h^{n(n-2) /4}\quad\hbox{and}\quad m_{L_{p}}(L_{p}\cap B^{\set{p}}_h)\ll_\epsilon h^{(n-2)^2 /4+\epsilon}
	$$
	for all $\epsilon>0$, and 
	$$
	\sup\left\{ \|\pi(\beta^{\set{p}}_h)\|:\, \pi \in \widehat G_{p}^{{\rm\tiny aut},0} \right\}	\gg_\epsilon h^{-(n-2)/2+\epsilon}\gg
	m_{G_{p}} (B^{\set{p}}_h)^{-2/n+\epsilon}.
	$$
	Similarly, when $n$ is odd, 
	$$
	m_{G_{p}}(B^{\set{p}}_h)\gg h^{(n-1)^2 /4}\quad\hbox{and}\quad m_{L_{p}}(L_{p}\cap B^{\set{p}}_h)\ll_\epsilon  h^{(n-1)(n-3) /4+\epsilon}
	$$
	for all $\epsilon>0$, and 
	$$
	\sup\left\{ \|\pi(\beta^{\set{p}}_h)\|:\, \pi \in \widehat G_{p}^{{\rm\tiny aut},0} \right\}	\gg_\epsilon h^{-(n-1)/2+\epsilon}\gg m_{G_{p}}(B^{\set{p}}_h)^{-2/(n-1)+\epsilon}.
	$$
	
	For ${\sf G}=\hbox{Sp}_{2n}$ with $n\ge 2$, 
	we apply the estimate \eqref{eq:est}
	with the subgroup ${\sf L}=\hbox{Sp}_{2(n-1)}$.
	We have
	$$
	m_{G_{p}}(B^{\set{p}}_h)\gg h^{n(n+1)}\quad\hbox{and}\quad m_{L_{p}}(L_{p}\cap B^{\set{p}}_h)\ll_\epsilon  h^{n(n-1) +\epsilon}
	$$
	for all $\epsilon>0$, and 
	$$
	\sup\left\{ \|\pi(\beta^{\set{p}}_h)\|:\, \pi \in \widehat G_{p}^{{\rm\tiny aut},0} \right\}	\gg_\epsilon h^{-2n +\epsilon} \gg m_{G_{p}}(B^{\set{p}}_h)^{-2/(n+1)+\epsilon}.
	$$

\end{proof}

\section{Lower bounds for operator norms in the discrete spectrum}\label{sec:descrete}

In this section we discuss lower bounds for the norms of averaging operators
for the discrete part of the automorphic spectrum and prove Theorem \ref{th:discrete}. We start by reviewing basic facts about the Hecke algebras and representations of $p$-adic groups (see, for instance, \cite{car}).
Let ${\sf G}$ be a classical simple simply connected algebraic $\Q$-group which is split over $\mathbb{Q}$. We fix a prime $p$ such that 
${\sf G}$ is isotropic and unramified over $\Q_p$.
Let $U_p$ be the hyperspecial maximal compact subgroup of $G_p$.
We consider the Hecke algebra $\mathcal{H}(G_p,U_p)$ 
consisting of compactly supported $U_p$-biinvariant functions on $G_p$ with the product defined by the convolution.
The structure of $\mathcal{H}(G_p,U_p)$ can be explicitly described as follows.
We fix a maximal split $\Q$-torus ${\sf T}$ of ${\sf G}$ 
and a Borel subgroup ${\sf B=TN}$ such that the Iwasawa decomposition 
$G_p=U_pT_pN_p$ holds.
Let $\Lambda_p:=T_p/(T_p\cap U_p)\simeq \Z^{\dim ({\sf T})}$
and $W:=N_{\sf G}({\sf T})/Z_{\sf G}({\sf T})$ be the Weyl group.
For $\beta\in \mathcal{H}(G_p,U_p)$, we set 
$$
\phi_{\beta}(t):=\Delta_p(t)^{1/2} \int_{N_p} \beta(tn)\, dm_{N_p}(n),\quad t\in \Lambda_p=T_p/(T_p\cap U_p),
$$
where $\Delta_p$ denotes the modular function of the group $T_pN_p$
and the invariant measure $m_{N_p}$ on $N_p$ is normalized so that 
$m_{N_p}(N_p\cap U_p)=1$. 
It is known that the map $\beta\mapsto \phi_{\beta}$
defines an algebra-isomorphism between the Hecke algebra $\mathcal{H}(G_p,U_p)$ and 
the  algebra $\C[\Lambda]^W$ of $W$-invariant polynomials in $\C[\Lambda]$. 
This allows to give a complete description of the spherical functions in terms of the unramified characters
$\chi:T_p\to\mathbb{C}^\times$. Given any such character, the corresponding 
spherical function is defined by
$$
\omega_\chi(g):=\int_{U_p} \Delta_p^{1/2}(t(gu))\chi(t(gu)) \, dm_{U_p}(u),\quad g\in G_p,
$$
where $t(\cdot)$ denotes the $T_p$-component with respect to the Iwasawa decomposition
$G_p=U_pT_pN_p$.
Moreover, two such spherical functions are equal iff the corresponding characters
are conjugate with respect to the Weyl group $W$.
We write $\mathcal{X}_p$ for the sets of unramified characters
$\chi:T_p\to\mathbb{C}^\times$ and $\mathcal{X}^{temp}_p$
for the subset of characters with $|\chi|=1$.
We additionally note that for this correspondence
$$
\chi(\phi_{\beta})=\beta(\omega_{\chi})\quad\hbox{for all $\chi\in \mathcal{X}_p$ and $\beta \in \mathcal{H}(G_p,U_p)$.}
$$
On the other hand, the spherical functions naturally arise 
from irreducible spherical unramified representations of the $G_p$. 
Given such a representation $\pi_p$, 
the corresponding spherical functions $\omega_\pi$ are defined by
$$
\omega_{\pi_p}(g):=\left<\pi_p(g)v_{\pi_p},v_{\pi_p}\right>,\quad g\in G_p,
$$
where $v_{\pi_p}$ is a unit-norm $U_p$-invariant vector, which is known to be unique up to scalar multiple. 
It follows from uniqueness that
$$
\pi_p(\beta)v_{\pi_p}=\beta(\omega_{\pi_p})v_{\pi_p}.
$$
Under the above identifications, the tempered
irreducible spherical unramified representations of $G_p$ correspond to the  characters in $\mathcal{X}^{temp}_p$, and
the Plancherel formula holds: for all $\beta \in\mathcal{H}(G_p,U_p)$,
$$
\int_{G_p} \abs{\beta(g)}^2\,dm_{G_p}(g)=
\int_{\chi\in \mathcal{X}_p^{temp}/W}\abs{\beta(\omega_\chi)}^2 d\nu_p(\chi),
$$
where $\nu_p$ denotes the normalised spherical Plancherel 
measure for the  group $G_p$ (see \cite[Th.2]{mac}).

Similarly, the spherical spectrum of $G_\infty$ is parametrized by a subset of 
$\frak{a}^*_{\mathbb{C}}$, where $\mathfrak{a}$ is the Lie algebra of the $\RR$-split torus in $G_\infty$, and $\frak{a}^*_{\mathbb{C}}$ denotes 
the complexificated dual space.
Then $i\frak{a}^*$ gives the parametrization of the tempered spherical spectrum.
Let $\Omega$ be a bounded domain in $i\mathfrak{a}^\ast$
with rectifiable boundary. The spherical Plancherel 
density associated to $\Omega$ is defined as
$$
\Lambda_{\Omega}(t):=C({G_\infty})\int_{t\Omega}\frac{{\bf c}(\rho)}{{\bf c}(\lambda)}d\lambda\,,$$
where $C({G_\infty})$ is an explicit positive constant,
$\bf c$ is the Harish-Chandra $\bf c$-function of $G_\infty$, and $\rho$ is half the sum of positive roots. 
It is known that 
\begin{equation}
\label{eq:aaas1}
\Lambda_{\Omega}(t)=C(\Omega) t^d +O(t^{d-1})
\end{equation}
with explicit $C(\Omega)>0$ and $d=d(G_\infty)\in\N$.

Let us now consider irreducible spherical unramified automorphic representations $\pi$ of $G_\infty\times G_p$
discretely embedded in $L_0^2(X_{p,1})$.
Such representations splits as a tensor product $\pi=\pi_{\infty}\otimes \pi_p$,
where  $\pi_{\infty}$ and $\pi_p$ denote 
the irreducible spherical representations
of $G_\infty$ and $G_p$ respectively.
We denote by 
$\lambda_{\pi_\infty}\in \frak{a}^*_{\mathbb{C}}$ the infinitesimal character of the Archimedean component $\pi_{\infty}$ and by 
$\lambda_{\pi_p}\in \mathcal{X}_p$ the characters corresponding to
the representations $\pi_p$ as described above.
The representation $L^2_0(X_{p,1})$ can be viewed as a subrepresenation of
$L^2_0({\sf G}(\mathbb{A})/{\sf G}(\Q))$ consisting of functions invariant under $\prod_{q\ne p} {\sf G}(\Z_p)$. 
Therefore, with this notation, the equdistribution result of \cite{mt,fm}
yields, in particular, that for all $\beta\in \mathcal{H}(G_p,U_p)$ and $t\ge 1$,
\begin{equation}\label{eq:equid}
\sum_{\pi :\lambda_{\pi_{\infty}}\in t\Omega } \chi_{\pi_p}(\phi_\beta) 
=\Lambda_\Omega(t) 
\int_{\chi\in \mathcal{X}_p^{temp}/W} \chi(\phi_{\beta})\, d\nu_p(\chi)
+O_\Omega \Big(\|\beta\|_{L^1}t^{d-\delta}\Big),
\end{equation}
with explicit $\delta>0$.
Here the sum is taken over irreducible discretely embedded spherical unramified automorphic representations.
Let us denote  the number of such representation 
with $\lambda_{\pi_{\infty}}\in t\Omega$ by $N(t)$.
Then taking $\beta=\chi_{U_p}$, we obtain from \eqref{eq:equid} that
\begin{equation}
\label{eq:aaas2}
N(t)=\Lambda_\Omega(t)+O_\Omega \big(t^{d-\delta}\big).
\end{equation}
The estimate \eqref{eq:equid} was proved in \cite{mt}
for $\hbox{SL}_d$ with somewhat weaker error term and in \cite{fm} for groups satisfying a technical condition, which is satisfied  in particular for all classical split groups. Formula \eqref{eq:equid} underlies the proof we give of 
Theorem \ref{th:discrete}, but we note that we will use only the existence of the limit, and not the effective error estimates formula \eqref{eq:equid} provides.

\begin{proof}[Proof of Theorem \ref{th:discrete}]
	We consider the representation
	$$
	\rho_t:=\bigoplus_{\pi: \lambda_{\pi_\infty}	\in t\Omega} \pi_p
	$$
	of $G_p$, where the sum is taken over spherical unramified representations $\pi=\pi_{\infty}\otimes \pi_p$ of $G_\infty\times G_p$
	discretely embedded in $L^2_0(X_{p,1})$. For each $\pi$, we denote by 
	$f_\pi$ the corresponding unit spherical vector and set
	$f_t:= \sum_{\pi: \lambda_{\pi_\infty}	\in t\Omega} f_\pi$.
	Then clearly
	\begin{equation}\label{upper-est}
	\norm{\rho_t(\beta)f_t}\le  \sum_{\pi:\lambda_{\pi_\infty}\in t\Omega } \norm{\rho_t(\beta)f_\pi}\le N(t)\norm{\rho_t(\beta)}\,.
	\end{equation}
	On the other hand, 
	$$\rho_t(\beta)f_\pi=\pi_{p}(\beta)f_\pi=\beta(\omega_{\pi_{p}})f_\pi.$$
	Since the different components in the sum defining $f_t$ are mutually orthogonal, we have
	\begin{align*}
	\norm{\rho_t(\beta)f_t}^2&=  \sum_{\pi:\lambda_{\pi_\infty}\in t\Omega } \norm{\pi_{p}(\beta)f_\pi}^2= \sum_{\pi:\lambda_{\pi_\infty}\in t\Omega } \abs{\beta(\omega_{\pi_{p}})}^2\\
	&= \sum_{\pi:\lambda_{\pi_\infty}\in t\Omega } (\beta^\ast\ast\beta)(\omega_{\pi_{p}})
	= \sum_{\pi:\lambda_{\pi_\infty}\in t\Omega } \chi_{\pi_p} (\phi_{\beta^\ast\ast\beta}).
	\end{align*}
	Therefore,
	\begin{align*}
	\norm{\rho_t(\beta)f_t}^2
	&=\Lambda_\Omega(t) 
	\int_{\chi\in \mathcal{X}_p^{temp}/W} \chi(\phi_{\beta^\ast\ast\beta})\, d\nu_p(\chi)
	+O_\Omega \Big(\|\beta^\ast\ast\beta\|_{L^1}t^{d-\delta}\Big)\\
	&=\Lambda_\Omega(t) 
	\int_{\chi\in \mathcal{X}_p^{temp}/W} |\beta(\omega_{\chi})|^2\, d\nu_p(\chi)
	+O_\Omega \Big(\|\beta\|^2_{L^1}t^{d-\delta}\Big).
	\end{align*}
	Hence, applying the Plancherel Formula, we deduce that
	$$
	\norm{\rho_t(\beta)f_t}^2
	=\Lambda_\Omega(t) \|\beta\|_{L^2}^2
	+O_\Omega \Big(\|\beta\|^2_{L^1}t^{d-\delta}\Big).
	$$
	Hence, using \eqref{eq:aaas1} and \eqref{eq:aaas2}, we conclude that 
	for all $\varepsilon > 0$ and sufficiently large $t$ (depending on $\beta$ and $\varepsilon$),
	$$
	\norm{\rho_t(\beta)f_t}\ge (1-\varepsilon) N(t)\,\|\beta\|_{L^2}.
	$$
	Comparing this estimate with \eqref{upper-est} we conclude that	
	$$
	\norm{\rho_{p,1}^{\text{cusp}}( \beta)}\ge \norm{\rho_t(\beta)}\ge (1-\varepsilon)\|\beta\|_{L^2},
	$$
	for all $\varepsilon> 0$. This completes the proof.
\end{proof}

\section{From discrepancy estimates to spectral gap}\label{sec:spectral bounds}

Let ${\sf G}$ be a simply connected $K$-simple linear algebraic group defined over a number field $K$ and $S$ a finite set of non-Archemedian completions of $K_v$. 
We use the notation introduced in Section \ref{sec:basic}, and 
in particular, we recall that $\Gamma_{S,\mathcal{O}}$ denote the family of 
congruence lattice subgroup in the product $G_\infty\times G_S$. 
When $S$ consists of unramified places, we denote by $U_S$ the hyperspecial maximal compact subgroup of $G_S$. 
Recall that $m_{G_S}$ and $m_{G_\infty}$ denote Haar measures on $G_\infty$ and $G_S$, respectively. We normalize the Haar measures so that $m_{G_S}(U_S)=1$ and $\Gamma_{S}$ has covolume one with respect to 
$m_{G_\infty}\times m_{G_S}$.

We will need the following lemma:

\begin{lem} \label{l:fund_dom}
Assume that $G_\infty$ is not compact.
Let $U$ be a compact open subgroup of $G_S$ and $\Gamma:=\Gamma_{S,\mathcal{O}}\cap U$. Let $F\subset G_\infty$ be a fundamental domain for $\Gamma$ in $G_\infty$.
Then $F\times U$ is a fundamental domain for $\Gamma_{S,\mathcal{O}}$ in $G_\infty\times G_S$.
\end{lem}

\begin{proof}
Since $G_\infty$ is not compact, and $G$ is assumed simply connected, it follows from the Strong Approximation Property
\cite[\S 7.4]{PlaRa} that the image of $\Gamma_{S,\cO}$ in $G_S$ is dense.
Consider an arbitrary $(g_\infty,g_S)\in G_\infty\times G_S$.
Then it follows from density that there exists $\gamma_1\in \Gamma_{S,\cO}$
such that $\gamma_1\in g_S^{-1}U_S$.
Hence, $(g_\infty\gamma_1 ,g_S\gamma_1)\in G_\infty\times U$.
Furthermore, since $F$ is a fundamental domain for 
$\Gamma$ in $G_\infty$, there exists $\gamma_2\in \Gamma$ such that 
$g_\infty\gamma_1\gamma_2\in F$. Then 
$(g_\infty\gamma_1 \gamma_2,g_S\gamma_1\gamma_2)\in F\times U$.
This proves that $G_\infty\times G_S=(F\times U)\Gamma_{S,\cO}$.

Suppose that $(F\times U)\gamma_1\cap (F\times U)\gamma_2\ne \emptyset$ for some $\gamma_1,\gamma_2\in \Gamma_{S,\cO}$.
Then it follows that $\gamma_1\gamma_2^{-1}\in \Gamma_{S,\cO}\cap U=\Gamma$,
and $F\gamma_1\gamma_2^{-1}\cap F\ne \emptyset$. Hence, we conclude that $\gamma_1=\gamma_2$. This completes the proof of the lemma.
\end{proof}

We consider the following counting problem for the congruence lattices $\Gamma_{S,\mathcal{O}}$.
Let us fix a left-invariant Riemannian metric on $G_\infty$.
For $x\in G_\infty$ and $r>0$, we denote by $B(x,r)$ the ball centered at $x$ of radius $r$ in $G_\infty$. Let $B$ be a bounded measurable subset of $G_S$
of positive measure.
We consider the counting function $\big|\Gamma_{S,\mathcal{O}}\cap (B(x,r)\times B)\big|$ and the discrepancy
$$
\mathcal{D}\big(\Gamma_{S,\mathcal{O}}, B(x,r)\times B\big):=\left| \frac{|\Gamma_{S,\mathcal{O}} \cap (B(x,r)\times B)|}{m_{G_S}(B)} - \frac{m_{G_\infty}(B(e,r))}{|\Gamma_S/\Gamma_{S,\cO}|} \right|.
$$
We shall show that estimates on the discrepancy can be used to establish  norm
bounds for $\|\pi(\beta)\|$ for the representations $\pi$ arising from the $G_S$-action on the spaces $X_{S,\cO}$, where $\beta$ is the probability measure supported on the set $B\subset G_S$.

We will be interested in utilizing $L^2$-bounds for the discrepancy function
as $x$ varies over $G_\infty$. 
Let us assume the subset $B$ is bi-invariant under
$U_S$. 
Then the discrepancy is left invariant
under the subgroup $\Gamma_{S,\mathcal{O}}^0:=\Gamma_{S,\mathcal{O}}\cap U_S$
(viewed as a subgroup of $G_\infty$), so that it defines a function on $\Gamma_{S,\mathcal{O}}^0\backslash G_\infty$. We set 
$$
E_{S,\mathcal{O}}(r,B):=\big\|\mathcal{D}\big(\Gamma_{S,\mathcal{O}}, B(\cdot ,r)\times B\big)\big\|_{L^2(\Gamma_{S,\mathcal{O}}^0\backslash G_\infty)},
$$
where $\Gamma_{S,\mathcal{O}}^0\backslash G_\infty$ is equipped with the invariant probability measure.

We now formulate and prove a more precise version of Theorem \ref{th:converse0}, as follows. 

\begin{thm}\label{th:converse}
Assume that $S$ consists of unramified places, and $G_\infty$ is not compact.
Let $B$ be a non-empty bounded $U_S$-bi-invariant subset of $G_S$.
Fix an irreducible unitary representation $\pi$ of $G_S$ which is discretely embedded in $\rho_{S,\mathcal{O}}$. Then there exists $r_0(\pi) > 0$ such that for $0 < r \le r_0(\pi)$,  the operator $\pi(\beta)$ satisfies the bound
$$\norm{\pi(\beta)}\le  2\abs{\Gamma_S/\Gamma_{S,\cO}} m_{G_\infty}(B(e,r))^{-1}E_{S,\mathcal{O}}(r,B).$$
\end{thm}

\begin{proof}
Consider  
$$
\chi_r(g_\infty,g_S):=\chi_{B(e,r)}(g_\infty)\chi_{U_S}(g_S),\quad\quad \hbox{for $(g_\infty,g_S)\in G_\infty\times G_S$},
$$
namely the characteristic function 
of the subset $B(e,r)\times U_S$ of $G_\infty\times G_S$.
Then 
$$
\phi_r(g_\infty, g_S):={\sum}_{\delta\in \Gamma_{S,\cO}} \chi_r(g_\infty\delta,g_S\delta)
$$ 
defines a function on the space $X_{S,\cO}=(G_\infty\times G_S)/\Gamma_{S,\cO}$.
By the Fubini--Tonelli Theorem, for $x\in G_\infty$ and $u\in U_S=U_S^{-1}$,
\begin{align*}
\int_{B}\phi_r\big(b^{-1}(x^{-1},u)\big)\, dm_{G_S}(b) &={\sum}_{\delta\in \Gamma_{S,\cO}} \int_{B} \chi_r(x^{-1}\delta, b^{-1}u\delta)\, dm_{G_S}(b)\\
& = {\sum}_{\delta\in \Gamma_{S,\cO}} \int_{B} \chi_{B(e,r)}(x^{-1}\delta)\chi_{U_S} (b^{-1}u\delta)  \, dm_{G_S}(b)\\
&=
{\sum}_{\delta\in \Gamma_{S,\cO} \cap (B(x,r)\times G_S)} m_{G_S}\big(u\delta U_S \cap B\big).
\end{align*}
Since the set $B$ is $U_S$-bi-invariant,
if $\delta\in B$, we have $u\delta U_S\subset B$, and
if $\delta\notin B$, we have $u\delta U_S\cap B=\emptyset$.
Hence, since $m_{G_S}(U_S)=1$, it follows that for every $u\in U_S$,
\begin{align}\label{eq:form}
\big|\Gamma_{S,\cO}\cap (B(x,r)\times B)\big|=\int_{B}\phi_r\big(b^{-1}(x^{-1},u)\big)\, dm_{G_S}(b).
\end{align}
We also compute:
\begin{align}
\int_{X_{S,\cO}} \phi_r \, d\mu_{S,\cO} &=
\int_{(G_\infty\times G_S)/\Gamma_{S,\cO}} \left(\sum_{\delta\in \Gamma_{S,\cO}} \chi_r(g\delta)\right)\, d\mu_{S,\cO}(g)\nonumber \\
&= \int_{G_\infty \times G_S} \chi_r(g)\, \frac{d(m_{G_\infty}\times m_{G_S})(g)}{|\Gamma_S/\Gamma_{S,\cO}|} \nonumber\\
&= \frac{m_{G_\infty}(B(e,r))m_{G_S}(U_S)}{|\Gamma_S/\Gamma_{S,\cO}|}  
=\frac{m_{G_\infty}(B(e,r))}{|\Gamma_S/\Gamma_{S,\cO}|}. \nonumber
\end{align}
Therefore, we conclude that for every $u\in U_S$,
\begin{align*}
\mathcal{D}(\Gamma_{S,\cO}, B(x,r)\times B) &=
\left| \frac{1}{m_{G_S}(B)}\int_{B}\phi_r\big(b^{-1}(x^{-1},u)\big)\, dm_{G_S}(b)-\int_{X_{S,\cO}} \phi_r\, d\mu_{S,\cO}  \right|\\
&=\left| \rho_{{S,\cO}}(\beta)\left( \phi_r -\int_{X_{S,\cO}} \phi_r\, d\mu_{S,\cO}\right)(x^{-1},u)  \right|.
\end{align*}

Let $\Omega$ be a measurable fundamental domain for
$\Gamma_{S,\mathcal{O}}^0:=\Gamma_{S,\cO}\cap U_S$  in $G_\infty$. 
Then for any $\Gamma_{S,\mathcal{O}}^0$-injective bounded measurable subset $Q$ of $G_\infty$
(in particular, for $Q=\Omega$), 
\begin{equation}
\label{eq:L2_0}
\left\| \rho_{{S,\cO}}(\beta)\left( \phi_r -\int_{X_{S,\cO}} \phi_r\, d\mu_{S,\cO}\right) \right\|_{L^2(Q^{-1}\times U_S)}
\le E_{S,\mathcal{O}}(r,B).
\end{equation}
By Lemma \ref{l:fund_dom}, $\Omega \times U_S$ is 
a fundamental domain for $\Gamma_{S,\cO}$ in $G_\infty\times G_S$.
Hence, we deduce from \eqref{eq:L2_0} that
\begin{equation}
\label{eq:L2_2}
\left\| \rho_{S,\mathcal{O}}(\beta) \phi_r -\int_{X_{S,\mathcal{O}}} \phi_r \, d\mu_{S,\mathcal{O}} \right\|_{L^2(X_{S,\mathcal{O}})}
\le E_{S,\mathcal{O}}(r,B).
\end{equation}
For $x\in G_\infty$ and $r>0$, let $\chi_{x,r}$ denote the characteristic function of the subset $B(x,r)\times U_S$ of $G_\infty\times G_S$.
Because of left-invariance of the metric, $\chi_{x,r}(g)=\chi_{r}(x^{-1}g)$.
We set 
$$
\phi_{x,\epsilon}(g):=\sum_{\gamma\in \Gamma_{S,\mathcal{O}}} \chi_{x,\epsilon}(g\gamma),
$$
which defines a function in $L^2(X_{S,\mathcal{O}})$.
For $x\in G_\infty$, let us consider the operators 
$$
\rho_{\infty,\mathcal{O}}(x):L^2(X_{S,\mathcal{O}})\to L^2(X_{S,\mathcal{O}}):\phi\mapsto\phi\circ x^{-1}.
$$
We observe that $\|\rho_{\infty,\mathcal{O}}(x)\|=1$, $\rho_{\infty,\mathcal{O}}(x)$ commutes with $\rho_{S,\mathcal{O}}(\beta)$,
and $\phi_{x,r}=\rho_{\infty,\mathcal{O}}(x)(\phi_r)$.
Hence, it follows from \eqref{eq:L2_2} that for every $x\in G_\infty$,
\begin{equation}
\label{eq:L2_3}
\left\| \rho_{S,\mathcal{O}}(\beta) \phi_{x,r} -\int_{X_{S,\mathcal{O}}} \phi_r \, d\mu_{S,\mathcal{O}} \right\|_{L^2(X_{S,\mathcal{O}})}
\le E_{S,\mathcal{O}}(r,B).
\end{equation}

We use this estimate to conclude the proof of the theorem as follows.
Let $\pi$ be an irreducible unitary representation of $G_S$ which is discretely embedded in $L_0^2(X_{S,\cO})$.
We observe that since $B$ is $U_S$-bi-invariant,
the image of $\pi(\beta)$ consists of $U_S$-invariant vectors.
Hence, if $\pi$ is not spherical, then $\pi(\beta)=0$.
Now suppose that $\pi$ is spherical and denote by $F_\pi$
the function  in $L_0^2(X_{S,\cO})$ which is the unique (up to a phase factor) 
$U_S$-invariant unit vector of $\pi$.
Then 
$$
\omega_\pi(g):=\left<\pi(g)F_\pi,F_\pi\right>,\quad\hbox{ with $g\in G_S$,}
$$
is
the spherical function associated to the representation $\pi$.
We have  
$$
\pi(\beta)F_\pi=\beta(\omega_\pi)F_\pi\quad\hbox{and}\quad \|\pi(\beta)\|^2=\|\pi(\beta*\beta^*)\|=\abs{\beta(\omega_\pi) }^2, 
$$
and every vector $F$ in the representation space of $\pi$ is an eigenvector of $\pi(\beta)$, namely $\pi(\beta)F=\lambda_F  F_\pi$. $F$ is orthogonal to the constant functions, and so it follows from \eqref{eq:L2_3} that when $F$ has unit norm 
\begin{align*}
\big|\left<\rho_{S,\mathcal{O}}(\beta^* * \beta)\phi_{x,r},F\right>\big|
&=\big|\left<\rho_{S,\mathcal{O}}(\beta)\phi_{x,r},\pi(\beta) F\right>\big|\\
&=\left|\left<\rho_{S,\mathcal{O}}(\beta)\phi_{x,r}-\int_{X_{S,\mathcal{O}}} \phi_r \, d\mu_{S,\mathcal{O}},\pi(\beta)F\right>\right|\\
&\le 
\left\| \rho_{S,\mathcal{O}}(\beta) \phi_{x,r} -\int_{X_{S,\mathcal{O}}} \phi_r \, d\mu_{S,\mathcal{O}} \right\|_{L^2(X_{S,\mathcal{O}})} \left\|\pi(\beta) F \right\|_{L^2(X_{S,\mathcal{O}})}\\
&\le  E_{S,\mathcal{O}}(r,B)\, \norm{\pi(\beta)}.
\end{align*}
On the other hand, let us choose a sequence of unit vectors $\psi_i$ in the 
representation $\pi$ such that 
$$
\|\pi(\beta^* * \beta)\psi_i\|\to \|\pi(\beta^* * \beta)\|. 
$$
Then $\pi(\beta^* * \beta)\psi_i=\lambda_i F_\pi$ with $\lambda_i\ge 0$ and $\lambda_i\to \|\pi(\beta^* * \beta)\|=\|\pi(\beta)\|^2$, and    
\begin{align*}
\big|\left<\rho_{S,\mathcal{O}}(\beta^* * \beta)\phi_{x,r},\psi_i\right>\big|=
\big|\left<\phi_{x,r}, \rho_{S,\mathcal{O}}(\beta^* *\beta)\psi_i\right>\big|=
\lambda_i \big|\left<\phi_{x,r}, F_\pi\right>\big|.
\end{align*}
Hence, we obtain the following norm bound
\begin{equation}
\label{eq:norm_low}
\|\pi(\beta)\|\le \big|\left<\phi_{x,r}, F_\pi\right>\big|^{-1} E_{S,\mathcal{O}}(r,B)
\end{equation}
provided that $\left<\phi_{x,r}, F_\pi\right>\ne 0$.

Let us now consider the function 
$$
f(g_\infty):=F_\pi\left(g_\infty\Gamma_{S,\mathcal{O}}\right)\quad\hbox{ with $g_\infty \in G_\infty$.}
$$
We note since $F_\pi$ is $U_S$-invariant,   we have $f(g_\infty)=F_\pi\left((g_\infty,u)\Gamma_{S,\mathcal{O}}\right)$
for all $u\in U_S$. 
In particular, $f$ is a well-defined measurable 
locally $L^2$-integrable function on $G_\infty$.
Since $G_\infty \Gamma_{S,\mathcal{O}}$ is dense in $G_\infty\times G_S$,
it is clear that $f\ne 0$.
We compute:
\begin{align*}
\left<\phi_{x,r}, F_\pi\right>&=\int_{X_{S,\mathcal{O}}} \left(\sum_{\gamma\in \Gamma_{S,\cO}} \chi_{x,r}(g\gamma)\right) \overline{F_\pi(g\Gamma_{S,\mathcal{O}})}\, d\mu_{S,\mathcal{O}}(g\Gamma_{S,\mathcal{O}})\\
&= \int_{G_\infty\times G_S}  \chi_{x,r}(g) \overline{F_\pi(g\Gamma_{S,\mathcal{O}})}\, \frac{d(m_{G_\infty}\times m_{G_S})(g)}{{|\Gamma_S/\Gamma_{S,\cO}|}}.
\end{align*}
We recall that $\chi_{x,r}$ is the characteristic function of the set $B(x,r)\times U_S$.
Since $F_\pi$ is $U_S$-invariant,
$$
\left<\phi_{x,r}, F_\pi\right>=|\Gamma_S/\Gamma_{S,\cO}|^{-1}\int_{B(x,r)} f\, dm_{G_\infty}.
$$
It follows from the Local Ergodic Theorem (see, for instance, \cite[Cor.~2.14]{mat}) that for almost every $x\in G_\infty$,
$$
\frac{1}{m_{G_\infty}(B(x,r))}\int_{B(x,r)} f\, dm_{G_\infty}\longrightarrow  f(x)\quad\hbox{ as $r\to 0^+$.}
$$
For a positive measure set of choices of $x$,
$$\abs{f(x)}\ge \frac 12 \norm{f}_\infty=\frac12\norm{F_\pi}_\infty\ge \frac12 \norm{F_\pi}_2=\frac 12\,,$$
  so that choosing a point $x$ where the foregoing inequality holds, and in addition convergence in the local ergodic theorem holds, it follows that for $0 < r < r_0(\pi) $ 
$$
\abs{\left<\phi_{x,r}, F_\pi\right>}\ge  \frac{1/2}{|\Gamma_S/\Gamma_{S,\cO}|} m_{G_\infty}(B(x,r)) \,.
$$
Therefore, the estimate \eqref{eq:norm_low} implies the claim of the theorem.
\end{proof}

\section{Spectral gap for forms of $\hbox{SL}_2$}\label{sec:sl2}

Let $\sf G$ be linear algebraic group defined over $\Q$ which is a form of  $\hbox{SL}_2$.
Namely, $\sf G$ can be realised as the group of norm one elements
of a quaternion algebra $\sf D$ defined over $\Q$:
$$
{\sf G}:=\{x\in {\sf D}:\, N(x)=1\}.
$$
Throughout this section, we always assume that $\sf G$ is isotropic over $\R$, 
or equivalently, ${\sf D}(\R)\simeq \hbox{M}_2(\R)$.
We fix an order $\Lambda$ of ${\sf D}(\Q)$ such that $N(\Lambda)\subset \mathbb{Z}$.
Then $N$ is an integral quadratic form with respect to this integral structure.
The $p$-adic norm $\|\cdot\|_p$ on ${\sf D}(\Q_p)$ is defined with respect to the order $\Lambda$.


The group $\Gamma:=G_\infty \cap \Lambda$ is an arithmetic lattice in 
$G_\infty$ corresponding to the integral structure defined by $\Lambda$. More generally, for a prime $p$, we consider the group 
$\Gamma_p:=G_\infty \cap \Lambda[p^{-1}]$ which is a lattice in the product $G_\infty\times G_p$. For $\ell\in\N$ coprime to $p$, we also consider 
the congruence subgroups $\Gamma_{p,\ell}:=\{\gamma\in \Gamma_p:\, \gamma=I\,(\hbox{mod}\,\ell) \}.$
The goal of this section is to analyze the unitary representations
$\rho_{p,\ell}$ of $G_p$ acting on the space $L_0^2(X_{p,\ell})$, where 
$X_{p,\ell}:=(G_\infty\times G_p)/\Gamma_{p,\ell}$. More specifically,
we will be interested in the averaging operators $\rho_{p,\ell}(\beta_h)$
defined with respect to the sets $B_h:=\{b\in G_p:\, \|b\|_p\le h\}.$

Let $m_p$ be the Haar measure on $G_p$ normalized so that $m_p({\sf G}(\Z_p))=1$,
and $m_\infty$ the Haar measure on $G_\infty$ normalized so that $m_\infty(G_\infty/\Gamma)=1$.

Our main result in this section is the following:

\begin{thm}\label{th:sl2}
Suppose that $\sf G$ is unramified 	over $\Q_p$.
Then for every $\ell\in \mathbb{N}$ coprime to $p$,
$$
\big\|\rho_{p,\ell}(\beta_h)\big\|\ll_{p,\ell} m_p(B_h)^{-\kappa},
$$
 where $\kappa=1/4$ if ${\sf G}$ is anisotropic over $\Q$,
and $\kappa=1/16$ if ${\sf G}$ is isotropic over $\Q$,
\end{thm}

It turns out that the representation-theoretic problem of estimating the norms
$\big\|\rho_{p,\ell}(\beta_h)\big\|$ is closely related to 
properties of the distribution of rational points contained in the quadratic surface $N(x)=1$.
For $h\in\N$ and a compactly supported function $w:{\sf D}(\mathbb{R})\to \mathbb{R}$,
we consider the counting function:
$$
\hbox{N}_h(N):=\sum_{x\in \Lambda: N(x)=h^2} w(h^{-1}x).
$$
More generally, for $\ell\in \N$ and a coset representative $\xi\in \Lambda/\ell\Lambda$, we define 
$$
\hbox{N}_h(N,w; \xi):=\sum_{x\in \xi+\ell\Lambda: N(x)=h^2} w(h^{-1}x).
$$
The behavior of $\hbox{N}_h(N,w; \xi)$ as $h\to\infty$ 
captures the distribution of the set of rational points $h^{-1}\Lambda$
with prescribed congruence condition.
In order to state an asymptotic formula for $N_h(N,w; \xi)$,
we need to introduce local densities.

For a compactly supported function $w$ on ${\sf D}(\RR)$, we define the {\it Archemedian local density} as
$$
\sigma_\infty(N,w):=\lim_{\epsilon\to 0^+} (2\epsilon)^{-1}\int_{|N(x)-1|\le \epsilon} w(x)\, dx,
$$
where the measure on $D(\RR)\simeq \R^4$ is normalized so that 
the lattice $\Lambda$ has covolume one.

Let $\ell=\prod_{q} {q}^{s_{q}}$ be the prime decomposition of $\ell$ and $\xi\in \Lambda/\ell\Lambda$.
We define, for a general integer $e\in \NN$,
$$
\hbox{N}_h(N, {q}^e, {q}^{s_{q}}, \xi):=\big|\{x\,\hbox{mod}\, {q}^{e+s_{q}}:\,\, x=\xi \,\hbox{mod}\, {q}^{s_{q}},\, N(x)=h^2\,\hbox{mod}\, {q}^e\}\big|.
$$
Then the {\it ${q}$-adic local density} is defined as
$$
\sigma_{q}(N,\xi,h):=\lim_{e\to\infty} \frac{N_h({q}^e,{q}^{s_{q}}, \xi)}{{q}^{3e}}.
$$
We also set
$$
\sigma_f(N,\xi,h):={\prod}_{q} \sigma_{q}(\xi,h).
$$

With these notations we state:

\begin{thm}\label{th:heath-brown}
For every smooth compactly supported function $w:{\sf D}(\mathbb{R})\to \mathbb{R}$
and $\xi \in \Lambda/\ell\Lambda$,
$$
\hbox{\rm N}_h(N,w; \xi)=\ell^{-4}\sigma_\infty(N,w)\sigma_f(N,\xi,h)h^2+O_{w,\ell,\epsilon}\big (h^{3/2+\epsilon}\big)\quad \hbox{ for all $\epsilon>0.$}
$$
\end{thm}

To prove Theorem \ref{th:heath-brown}, we follow closely the refined circle method for quadratic forms developed by Heath-Brown \cite{HB96}, for the case of forms in four variables. This  method is particularly suitable for the derivation of estimates which are uniform over families of functions $w$, which will be crucial for our results, cf. Theorem \ref{th:heath-brown2} below. We note however that in order to establish uniform bounds on $\big\|\rho_{p,\ell}(\beta_h)\big\|$ over $\ell$, we have to keep track of additional congruence conditions beyond those considered in \cite{HB96}. This causes considerable technical complications, which motivates our decision to give a full account of the necessary arguments below.  

Making the identification $\Lambda\simeq \mathbb{Z}^4$,
we view the reduced norm $N$ as an integral quadratic form in four variable.
We denote by $|\cdot|$ the Euclidean norm on ${\sf D}(\R)\simeq \R^4
$ defined by this identification.

The starting form of our argument is a convenient representation of  
the Dirac function  $\delta_n=1$ if $n=0$ and $\delta_n=0$ for $n\ne 0$.
We recall (see \cite[Th.~1]{HB96} and \cite{dfi}) that for every $Q>1$,
\begin{equation}
\label{eq:d_n}
\delta_n=c_Q Q^{-2} \sum_{k=1}^\infty {\sum_{a\,\hbox{\tiny mod}\, k}}^{\!\!\!*} e_k(an) H(Q^{-1}k,Q^{-2}n),
\end{equation}
where $c_Q=1+O_N(Q^{-N})$ for any $N>0$, 
the sum is taken over $a$ coprime to $k$,
$e_k(x)=\exp(2\pi ix/k)$, 
and $H$ is a certain explicit smooth function on $(0,\infty)\times \mathbb{R}$.
To simplify notation, we set $F_h(x)=N(x)-h^2$.
It follows from \eqref{eq:d_n} that
$$
\hbox{N}_h(N,w; \xi)=c_h h^{-2} \sum_{x\in \xi+\ell\mathbb{Z}^4}\sum_{k=1}^\infty {\sum_{a\,\hbox{\tiny mod}\, k}}^{\!\!\!*} w(h^{-1}x)e_k(a F_h(x)) H\big(h^{-1}k,h^{-2}F_h(x)\big),
$$
and we are required to estimate the following sum
\begin{align*}
&\sum_{x\in \xi+\ell\mathbb{Z}^4} w(h^{-1}x)e_k(a F_h(x)) H\big(h^{-1}k,h^{-2}F_h(x)\big)\\
=&\sum_{{\tiny b\,\hbox{\tiny mod}\, k}} \sum_{{\tiny \begin{tabular}{c} $z\,\hbox{\tiny mod}\, k\ell$\\ $z=\xi(\ell), z=b(k)$\end{tabular}}}  e_k(a F_h(b)) 
\left(\sum_{y\in \mathbb{Z}^4} f(y) \right),
\end{align*}
where 
$$
f(y):=w(h^{-1}(z+(k\ell)y))H\big(h^{-1}k,h^{-2}F_h(z+(k\ell)y)\big)
$$
is a smooth compactly supported function.
By the Poisson Summation Formula,
$$
\sum_{y\in \mathbb{Z}^4} f(y)=\sum_{c\in \mathbb{Z}^4} \hat f(c),
$$
where
$$
\hat f(c)=\int_{\mathbb{R}^4} f(y) e^{-2\pi i (c\cdot y)}\,dy
=(k\ell)^{-4} e_{k\ell}(c\cdot z) I_{k,\ell}(c)
$$
with
\begin{align*}
I_{k,\ell}(c)&:=\int_{\mathbb{R}^4} w(h^{-1}x) H\big(h^{-1}k, h^{-2}F_h(x)\big)e_{k\ell}(-c\cdot x)\, dx\\
&=h^4\int_{\mathbb{R}^4} w(x) H\big(h^{-1}k, N(x)-1\big)e_{k\ell}(-hc\cdot x)\, dx.
\end{align*}
Hence, we deduce that
\begin{equation}\label{eq:N_h_2}
\hbox{N}_h(N,w; \xi)=c_h h^{-2} \sum_{c\in\mathbb{Z}^4}\sum_{k=1}^\infty 
(k\ell)^{-4} S_{k}(c;\xi) I_{k,\ell}(c),
\end{equation}
where
\begin{align*}
S_{k}(c;\xi):=&{\sum_{a\,\hbox{\tiny mod}\, k}}^{\!\!\!*} \sum_{{\tiny b\,\hbox{\tiny mod}\, k}}
\sum_{{\tiny \begin{tabular}{c} $z\,\hbox{\tiny mod}\, k\ell$\\ $z=\xi(\ell), z=b(k)$\end{tabular}}}  e_{k\ell}(a\ell F_h(b) +c\cdot z)\\
=&{\sum_{a\,\hbox{\tiny mod}\, k}}^{\!\!\!*}
\sum_{{\tiny \begin{tabular}{c} $z\,\hbox{\tiny mod}\, k\ell$\\ $z=\xi(\ell)$\end{tabular}}} e_{k\ell}(a\ell F_h(z) +c\cdot z)\\
=&{\sum_{a\,\hbox{\tiny mod}\, k}}^{\!\!\!*}e_{k}(-a h^2)
\sum_{{\tiny \begin{tabular}{c} $z\,\hbox{\tiny mod}\, k\ell$\\ $z=\xi(\ell)$\end{tabular}}} e_{k\ell}(a\ell N(z) +c\cdot z).
\end{align*}

First, we observe that the sum $S_{k}(c;\xi)$ has the following multiplicative property:
\begin{lem} \label{l:mult}
Let $k=k_1k_2$ and $\ell=\ell_1\ell_2$ such that $k_1\ell_1$ is coprime to $k_2\ell_2$.
Choose integers $\bar k_1,\bar k_2,\bar \ell_1,\bar\ell_2$
such that 
$$
k_1\bar k_1=1\,(k_2\ell_2),\;\; k_2\bar k_2=1\,(k_1\ell_1),\;\;
\ell_1\bar \ell_1=1\,(k_2\ell_2),\;\; \ell_2\bar \ell_2=1\,(k_1\ell_1).
$$
Then
$$
S_k(c;\xi\,(\ell))=S_{k_1}(\bar k_2\bar \ell_2 c;\xi\,(\ell_1))\, S_{k_2}(\bar k_1\bar \ell_1 c;\xi\,(\ell_2)).
$$
\end{lem}

\begin{proof}
We also choose integers $\bar \ell_1$ and $\bar \ell_2$ such that
$\ell_1\bar \ell_1=1\,(k_2\ell_2)$ and $k_2\bar k_2=1\,(k_1\ell_1)$.
We write 
$$
z=k_2\ell_2\bar k_2\bar \ell_2 z_1+k_1\ell_1\bar k_1\bar \ell_1 z_2.
$$
If $z_1$ runs over the integral vectors modulo $k_1\ell_1$ such that
$z_1=\xi\,(\ell_1)$, and $z_2$ runs over the integral vectors modulo $k_2\ell_2$ such that
$z_2=\xi\,(\ell_2)$, then $z$ runs precisely over 
the integral vectors modulo $k\ell$ such that $z=\xi\,(\ell)$.
Then since the map $N$ is a quadratic form, 
$$
N(k_2\ell_2\bar k_2\bar \ell_2 z_1+k_1\ell_1\bar k_1\bar \ell_1 z_2)=
(k_2\ell_2\bar k_2\bar \ell_2)^2 N(z_1)+ (k_1\ell_1\bar k_1\bar \ell_1)^2 N(z_2)\;(\hbox{mod}\; k_1k_2\ell),
$$
and
\begin{align*}
e_{k\ell}(a\ell N(z))
&=e_{k_1k_2\ell }\big(a(\ell_1\ell_2(k_2\ell_2\bar k_2\bar \ell_2)^2 N(z_1)+ \ell_1\ell_2(k_1\ell_1\bar k_1\bar \ell_1)^2 N(z_2))\big)\\
&=e_{k_1\ell_1}\big(a \ell_1 k_2\ell_2^2\bar k_2^2 \bar\ell_2^2 N(z_1)\big)e_{k_2\ell_2}\big(a \ell_2 k_1\ell_1^2\bar k_1^2 \bar\ell_1^2 N(z_2)\big)\\
&=e_{k_1\ell_1}\big(a  \bar k_2 \ell_1 N(z_1)\big)e_{k_2\ell_2}\big(a  \bar k_1\ell_2 N(z_2)\big).
\end{align*}
Similarly,
\begin{align*}
e_{k\ell}(c\cdot z)&=e_{k_1\ell_1k_2\ell_2}\big(c\cdot k_2\ell_2\bar k_2\bar \ell_2 z_1+c\cdot k_1\ell_1\bar k_1\bar \ell_1 z_2\big)\\
&=e_{k_1\ell_1}\big(( \bar k_2\bar \ell_2 c)\cdot  z_1\big) 
e_{k_2\ell_2}\big((\bar k_1\bar \ell_1 c)\cdot z_2\big).
\end{align*}
Setting
\begin{align*}
S_{k}(a,k,\xi\,(\ell)):=\sum_{{\tiny \begin{tabular}{c} $z\,\hbox{\tiny mod}\, k\ell$\\ $z=\xi(\ell)$\end{tabular}}} e_{k\ell}(a\ell N(z) +c\cdot z),
\end{align*}
we conclude that
$$
S_{k}(a,k,\xi\,(\ell))=
S_{k_1}\big(a\bar k_2, \bar k_2\bar \ell_2 c;\xi\,(\ell_1)\big)\,
S_{k_2}\big(a \bar k_1, \bar k_1\bar \ell_1 c;\xi\,(\ell_2)\big).
$$
Then
\begin{align*}
S_k(c;\xi\,(\ell))
={\sum_{a\,\hbox{\tiny mod}\, k}}^{\!\!\!*} e_{k}(-a h^2)
S_{k_1}\big(a\bar k_2, \bar k_2\bar \ell_2 c;\xi\,(\ell_1)\big)\,
S_{k_2}\big(a\bar k_1, \bar k_1\bar \ell_1 c;\xi\,(\ell_2)\big).
\end{align*}
We observe that every residue $a \,\hbox{mod}\, k$
coprime to $k$ can be uniquely written as $k_2a_1+k_1a_2$,
where $a_i$ is a residue modulo $k_i$ coprime to $k_i$.
Hence, the above sum can be rewritten as:
\begin{align*}
&{\sum_{a_1\,\hbox{\tiny mod}\, k_1}}^{\!\!\!\!*}\;{\sum_{a_2\,\hbox{\tiny mod}\, k_2}}^{\!\!\!\!*} e_{k}(- (k_2a_1+k_1a_2) h^2)
S_{k_1}\big((k_2a_1+k_1a_2)\bar k_2, \bar k_2\bar \ell_2 c;\xi\,(\ell_1)\big)\\
&\quad\quad\quad\quad\quad\quad\quad  \times 
S_{k_2}\big( (k_2a_1+k_1a_2)\bar k_1, \bar k_1\bar \ell_1 c;\xi\,(\ell_2)\big)\\
=&{\sum_{a_1\,\hbox{\tiny mod}\, k_1}}^{\!\!\!\!*}\;{\sum_{a_2\,\hbox{\tiny mod}\, k_2}}^{\!\!\!\!*} e_{k_1}(- a_1 h^2) e_{k_2}(- a_2 h^2)
S_{k_1}\big(a_1, \bar k_2\bar \ell_2 c;\xi\,(\ell_1)\big)
S_{k_2}\big(a_2, \bar k_1\bar \ell_1 c;\xi\,(\ell_2)\big)\\
=&\;S_{k_1}(\bar k_2\bar \ell_2 c;\xi\,(\ell_1))\,S_{k_2}(\bar k_1\bar \ell_1 c;\xi\,(\ell_2)).
\end{align*}
This proves the lemma.
\end{proof}

We record the required properties of the term $S_k(c;\xi)$:

\begin{lem}[properties of the sum $S_k\big(c;\xi\big)$]
	\label{l:sq}
For every $c\in \mathbb{Z}^4$ and $\xi\in \Lambda/\ell\Lambda$,	
\begin{enumerate}
	\item[(i)] $|S_k\big(c;\xi\big)|\ll_{\ell} k^{3}$,
	\item[(ii)] For all $\epsilon>0$,
	\begin{align*}
	\sum_{k\le X} |S_k(c;\xi)| &\ll_{\ell,\epsilon} X^{7/2+\epsilon} |c|^\epsilon\quad\hbox{when $c\ne 0$,}
	\end{align*}
	and
	\begin{align*}
	\sum_{k\le X} |S_k(0;\xi)| &\ll_{\ell,\epsilon} X^{7/2+\epsilon} h^\epsilon.
	\end{align*}
\end{enumerate}
\end{lem}

\begin{proof}
The proof of (i) proceeds as \cite[Lem.~25]{HB96} with minor modifications.
Applying the Cauchy--Schwarz inequality to the sum
$$
S_k\big(c;\xi\big)={\sum_{a\,\hbox{\tiny mod}\, k}}^{\!\!\!*}e_{k}(-a h^2)
\sum_{{\tiny \begin{tabular}{c} $z\,\hbox{\tiny mod}\, k\ell$\\ $z=\xi(\ell)$\end{tabular}}} e_{k\ell}(a\ell N(z) +c\cdot z),
$$
we obtain that 
\begin{align}\label{eq:cancel}
|S_k\big(c;\xi\big)|^2&\le \phi(k)
{\sum_{a\,\hbox{\tiny mod}\, k}}^{\!\!\!*}
\left|\sum_{{\tiny \begin{tabular}{c} $z\,\hbox{\tiny mod}\, k\ell$\\ $z=\xi(\ell)$\end{tabular}}} e_{k\ell}\Big(a\ell N(z) +c\cdot z\Big)\right|^2\\
&\le \phi(k)
{\sum_{a\,\hbox{\tiny mod}\, k}}^{\!\!\!*}
\sum_{{\tiny \begin{tabular}{c} $z_1,z_2\,\hbox{\tiny mod}\, k\ell$\\ $z_1=z_2=\xi(\ell)$\end{tabular}}} e_{k\ell}\Big(a\ell (N(z_2)-N(z_1)) +c\cdot (z_2-z_1)\Big).\nonumber
\end{align}
We write $N(z)={}^t z A z$ for a symmetric matrix $A$ and $z_2=z_1+v$. Then 
the last sum can be rewritten as
$$
\sum_{{\tiny \begin{tabular}{c} $z_1,v\,\hbox{\tiny mod}\, k\ell$\\ $z_1=\xi(\ell)$, $v=0(\ell)$\end{tabular}}} e_{k\ell}\big(a\ell N(v) +c\cdot v\big)e_{k}\big(2a\,{}^t z_1 A v\big).
$$
We consider first the sum over the residues $z_1$ such that the vector 
$2\ell \,{}^t z_1 A$ satisfies $2\ell\,{}^t z_1 A=0\,(\hbox{mod}\, k)$.
It is clear that the number of such residues $z_1\,\hbox{mod}\, k\ell $ 
with $z_1=\xi(\ell)$ is $O_{A,\ell}(1)$.
Since the sum over $v$ has $k^4$ terms,
this implies that the contribution of the sum over such $z_1$
to \eqref{eq:cancel} is at most $O_{A,\ell}(k^6)$.
Now we consider the sum over the residues $z_1$
such that $2\ell\,{}^t z_1 A\ne 0\,(\hbox{mod}\, k)$.
This implies that 
$$
\sum_{{\tiny \begin{tabular}{c} $v\,\hbox{\tiny mod}\, k\ell$\\  $v=0(\ell)$\end{tabular}}} e_{k}\big(2a\,{}^t z_1 A v\big)=
\sum_{v'\,\hbox{\tiny mod}\, k} e_{k}\big(2a\ell \,{}^t z_1 A v'\big)=0.
$$
Therefore, we conclude that 
$$
|S_k\big(c;\xi\big)|=O_{A,\ell}(k^3),
$$
which proves (i).
	
To prove (ii), let us first consider the case when $k$ is coprime to $\ell$. Then
\begin{align*}
S_{k}(c;\xi)
=&{\sum_{a\,\hbox{\tiny mod}\, k}}^{\!\!\!*}
\sum_{{\tiny \begin{tabular}{c} $z\,\hbox{\tiny mod}\, k\ell$\\ $z=\xi(\ell)$\end{tabular}}} e_{k\ell}\big(a\ell F_h(z) +c\cdot z\big)\\
=&{\sum_{a\,\hbox{\tiny mod}\, k}}^{\!\!\!*}\;
\sum_{z'\,\hbox{\tiny mod}\, k} e_{k\ell}\big(a\ell F_h(\xi+\ell z') +c\cdot (\xi+\ell z')\big).
\end{align*}
Let $\bar \ell$ be the residue modulo $k$ such that $\ell\bar \ell=1\,(\hbox{mod}\, k)$. Then it follows that
\begin{align*}
S_{k}(c;\xi)
=&{\sum_{a\,\hbox{\tiny mod}\, k}}^{\!\!\!*}\;
\sum_{z''\,\hbox{\tiny mod}\, k} e_{k\ell}\big(a\ell F_h(z'') +c\cdot (\xi+\ell \bar \ell (z''-\xi))\big),
\end{align*}
and 
\begin{equation}\label{eq:b1}
|S_{k}(c;\xi)|\le |S_{k}(\bar \ell c)|,
\end{equation}
where
$$
S_{k}(c):={\sum_{a\,\hbox{\tiny mod}\, k}}^{\!\!\!*}
\sum_{z\,\hbox{\tiny mod}\, k} e_{k}\big(a F_h(z) +c\cdot z\big).
$$
The quantities $S_{k}(c)$ were investigated in \cite{HB96}.

We proceed as in \cite[Lem.~28]{HB96}.
Let us decompose $k=k_1k_2$ where $k_1$ is square-free, $k_2$ is square-full and $(k_1,k_2)=1$. We also decompose $\ell=\ell_1\ell_2$ so that $k_1\ell_1$ is coprime to 
$k_2\ell_2$. Then according to Lemma \ref{l:mult},
$$
S_k(c;\xi\,(\ell))=S_{k_1}(\bar k_2\bar \ell_2 c;\xi\,(\ell_1))\, S_{k_2}(\bar k_1\bar \ell_1 c;\xi\,(\ell_2)).
$$
We estimate the second factor using (i):
$$
S_{k_2}(\bar k_1\bar \ell_1 c;\xi\,(\ell_2))\ll_{\ell} k_2^3.
$$
The first factor can be decomposed further into product over primes $q|k_1$.
We use that according to Lemma \cite[Lem.~26]{HB96}, for all primes $q$ which do not divide the discriminant of $N$ and do not divide both $h^2$ and $A^{-1}c$, 
\begin{align}\label{l:b2}
|S_q(c)|\ll q^{5/2} (q,h^2, A^{-1}c )^{1/2}.
\end{align}
We note that this bound also holds when $q$ divide both $h^2$ and $A^{-1}c$
because of the estimate (i). Now for finitely many primes $q$ dividing either $q$ or $\ell$, we use the bound from (i), and for the remaining primes we use the estimate \eqref{l:b2} taking \eqref{eq:b1} into account. This gives the estimate:
\begin{align*}
\sum_{k\le X} |S_k(c;\xi)| &\le {\sum_{k_2\le X}}^{\!\!\!'}\; {\sum_{k_1\le X/k_2}}^{\!\!\!\!\!''} \big|S_{k_1}(\bar k_2\bar \ell_2 c;\xi\,(\ell_1))\big|\, \big|S_{k_2}(\bar k_1\bar \ell_1 c;\xi\,(\ell_2))\big|\\
&\ll_\ell {\sum_{k_2\le X}}^{\!\!\!'}\; {\sum_{k_1\le X/k_2}}^{\!\!\!\!\!''} 
k_1^{5/2} (k_1,h^2, A^{-1}c )^{1/2}\, k_2^3 \\
&\le X^{5/2 }{\sum_{k_2\le X}}^{\!\!\!'}\; k_2^{1/2 }{\sum_{k_1\le X/k_2}}^{\!\!\!\!\!''} 
(k_1,h^2, A^{-1}c )^{1/2}. 
\end{align*}
Here the sums are taken over square-full $k_2$ and square-free $k_1$ respectively.
We use that for $k \ne 0$, 
$$
\sum_{n\le Y} (n,k)\le Y d(k)\ll_{\epsilon} Y k^\epsilon
$$
for all $\epsilon>0$. Hence, we conclude that if $c\ne 0$,
\begin{align*}
\sum_{k\le X} |S_k(c;\xi)| \ll_{\ell,\epsilon} X^{7/2 }|c|^\epsilon {\sum_{k_2\le X}}^{\!\!\!'}\; k_2^{-1/2 } 
\ll_{\ell,\epsilon} X^{7/2+\epsilon }|c|^\epsilon.
\end{align*}
Similarly, we obtain 
\begin{align*}
\sum_{k\le X} |S_k(0;\xi)| \ll_{\ell,\epsilon} X^{7/2+\epsilon }h^\epsilon.
\end{align*}
This completes the proof of (ii).
\end{proof}

We shall use the following properties of the integral $I_{k,\ell}(c)$:

\begin{lem}[properties of the integral $I_{k,\ell}(c)$]
	\label{l:I_q}
	
	\begin{enumerate}
	\item[(i)] $I_{k,\ell}(c)=0$ for all $k\ge h$.
	\item[(ii)] For every $c\ne 0$ and $k>0$,
	$$
	I_{k,\ell}(c)\ll_{w,\ell, k} h^5 k^{-1}|c|^{-k}.
	$$
	\item[(iii)] For $c\ne 0$ and $\epsilon\in (0,1/2)$,
	$$
	I_{k,\ell}(c)\ll_{w,\ell,\epsilon} \frac{h^{3+2\epsilon}k^{1-2\epsilon}}{|c|^{1-\epsilon}}.
	$$
	\item[(iv)] $I_{k,\ell}(0)\ll_{w,\ell} h^4$.
	
	\item[(v)] For every $k\ll h$ and $k>0$,
	$$
	I_{k,\ell}(0)=h^4\Big( \sigma_\infty(N,w)+O_{w,\ell,k}\big((k/h)^k\big)\Big).
	$$
	\end{enumerate}
\end{lem}

\begin{proof}
Part (i) follows from properties of the function $H$ (see \cite[Lem.~4]{HB96}).
We observe that
\begin{align*}
I_{k,\ell}(c)
=(h\ell)^4\int_{\mathbb{R}^4} w(\ell x) H\big(h^{-1}k, \ell^2 N(x)-1\big)e_{k}(-hc\cdot x)\, dx.
\end{align*}
The properties of the integrals
\begin{align*}
I_{k}(c):=h^4\int_{\mathbb{R}^4} w(x) H\big(h^{-1}k, N(x)-1\big)e_{k}(-hc\cdot x)\, dx.
\end{align*}
have been studied in \cite{HB96}, and properties (ii)--(iv) can be deduced directly from there.
In particular, (ii) follows from \cite[Lem.~19]{HB96}, (iii) from \cite[Lem.~22]{HB96}, and (iv) from \cite[Lem.~22]{HB96}.
To prove (v), we use \cite[Lem.~13]{HB96}.
Here we also use that $\sigma_\infty(\ell^2 N,w\circ \ell)= \ell^{-4}\sigma_\infty( N,w)$.
\end{proof}

\begin{proof}[Proof of Theorem \ref{th:heath-brown}]
Our starting point is the formula
$$
\hbox{N}_h(N,w; \xi)=c_h h^{-2} \sum_{c\in\mathbb{Z}^4}\sum_{k=1}^\infty 
(k\ell)^{-4} S_{k}(c;\xi) I_{k,\ell}(c).
$$	
Using Lemmas \ref{l:sq}(i) and  \ref{l:I_q}(ii), we deduce that
for every $\theta>0$,
$$
\left|\sum_{|c|>h^\epsilon}\sum_{k=1}^\infty
(k\ell)^{-4} S_{k}(c;\xi) I_{k,\ell}(c)\right| \ll_{w,\ell,k} h^5 \left(\sum_{k=1}^\infty k^{-2}\right)\left(\sum_{|c|>h^\epsilon} |c|^{-\theta}\right)
=O_{w,\ell,\theta}(1),
$$
when $k$ is chosen to be sufficiently large. 
Now it is sufficient to consider the sum over $|c|\le h^\epsilon$.
It follows from Lemma \ref{l:I_q}(iii) that when $c\ne 0$ with $|c|\le h^\epsilon$,
$$
\left|\sum_{R\le k< 2R} (k\ell)^{-4} S_k(c;\xi)I_{k,\ell}(c)\right| 
\ll_{w,\ell,\epsilon} h^{3+2\epsilon}  \sum_{R\le k< 2R} k^{-3} |S_k(c;\xi)|,
$$
Here we used that $|c|\ge 1$ for $c\in\mathbb{Z}^4\{0\}$.
Hence, using Lemma \ref{l:sq}(ii) and summation by parts, we deduce that
$$
\sum_{R\le k< 2R} k^{-3} |S_k(c;\xi)| \ll_{\ell,\epsilon} R^{1/2+2\epsilon}.
$$
Hence, we deduce that for every $c\ne 0$ with $|c|\le h^\epsilon$,
$$
\left|\sum_{k=1}^h (k\ell)^{-4} S_k(c;\xi)I_{k,\ell}(c)\right| 
\ll_{w,\ell,\epsilon} h^{7/2+4\epsilon}.
$$
Moreover, according to Lemma \ref{l:I_q}(i), 
only the terms with $k\le h$ are non-zero.
Therefore, we conclude that for every $\epsilon>0$,
\begin{equation}
\label{eq:N_h2}
\hbox{N}_h(N,w; \xi)=c_h h^{-2}\ell^{-4} \sum_{k=1}^\infty 
k^{-4} S_{k}(0;\xi) I_{k,\ell}(0)+
O_{w,\ell,\epsilon} \big(h^{3/2+\epsilon}\big).
\end{equation}
Using Lemma \ref{l:sq}(iv), Lemma \ref{l:I_q}(ii), and summation by parts,
we deduce that for every $\epsilon>0$,
$$
\sum_{R<k\le 2R} 
k^{-4} S_{k}(0;\xi) I_{k,\ell}(0)\ll_{w,\ell,\epsilon} h^{4+\epsilon} R^{-1/2+2\epsilon}.
$$
This implies that for every $\epsilon>0$,
$$
\sum_{k>h^{1-\epsilon}} 
k^{-4} S_{k}(0;\xi) I_{k,\ell}(0)\ll_{w,\ell,\epsilon} h^{7/2+\epsilon}.
$$
It follows from Lemma \ref{l:I_q}(v) that
$$
\sum_{k\le h^{1-\epsilon}} 
k^{-4} S_{k}(0;\xi) I_{k,\ell}(0)
=h^4\sigma_\infty(N,w) \left(\sum_{k\le h^{1-\epsilon}} 
k^{-4} S_{k}(0;\xi) \right)+O_{w,\ell,\epsilon}(1).
$$
Hence, we deduce from \eqref{eq:N_h2} that
$$
\hbox{N}_h(N,w; \xi)=h^{2}\sigma_\infty(N,w)\ell^{-4} \left( \sum_{k\le h^{1-\epsilon}} 
k^{-4} S_{k}(0;\xi)\right) +
O_{w,\ell,\epsilon} \big(h^{3/2+\epsilon}\big).
$$
The last sum corresponds to the classical singular series. It follows from Lemma \ref{l:sq}
that it converges absolutely, and moreover, 
$$
\sum_{k\le X} 
k^{-4} S_{k}(0;\xi)=\sum_{k=1}^\infty 
k^{-4} S_{k}(0;\xi)+O_{\ell,\epsilon}\big(X^{-1/2+\epsilon}\big)
$$
for all $\epsilon>0$.
It is well-known that if the singular series converges,
it is equal to the product of local densities.
\end{proof}

We will also need a uniform version of Theorem \ref{th:heath-brown} which
applies to families of functions $w_g(x):=w(g^{-1}x)$ with $g\in G_\infty$.

\begin{thm}\label{th:heath-brown2}
	Let $w:D_\infty\to \mathbb{R}$ be a smooth compactly supported function,
	$\xi \in \Lambda/\ell\Lambda$, and $g\in G_\infty$.
	Then for every $\delta>0$, $\theta>4$, and $\epsilon>0$,
	$$
	\hbox{\rm N}_h(N,w_g; \xi)=\ell^{-4}\sigma_\infty(N,w)\sigma_f(N,\xi)h^2
	+O_{w,\ell,\delta,\theta,\epsilon}\Big(\|g\|^\theta h^{3-(\theta-4)\delta}+\|g\| h^{3/2+3\delta+\epsilon} \Big).
	$$
\end{thm}

\begin{proof}
As in the proof of Theorem \ref{th:heath-brown},
$$
\hbox{N}_h(N,w_g; \xi)=c_h h^{-2} \sum_{c\in\mathbb{Z}^4}\sum_{k=1}^\infty 
(k\ell)^{-4} S_{k}(c;\xi) I_{g,k,\ell}(c),
$$
where
\begin{align*}
I_{g,k,\ell}(c)
:=h^4\int_{\mathbb{R}^4} w(g^{-1}x) H(h^{-1}k, N(x)-1)e_{k\ell}(-hc\cdot x)\, dx.
\end{align*}
Since the norm $N$ is $G_\infty$-invariant, we obtain
\begin{align*}
I_{g,k,\ell}(c)
=h^4\int_{\mathbb{R}^4} w(x) H(h^{-1}k, F(x))e_{k\ell}(-hc\cdot gx)\, dx
=I_{k,\ell}({}^tg c).
\end{align*}
Our argument proceeds as in  the proof of Theorem \ref{th:heath-brown}
taking the dependences on $g$ into account. Throughout the proof,
we will have to deal with the maps ${\sf D}(\R)\to {\sf D}(\R):x\mapsto 
{}^tg^{\pm 1} x$ with $g\in G_\infty$. We use the norm of these maps are estimated
in terms of the Euclidean norm $\|g\|_E:=\big(\sum_{i,j} g_{ij}^2 \big)^{1/2}$ on the group $G_\infty\simeq \hbox{SL}_2(\R)$.

First, we consider the terms with $|c|>h^\delta$. By Lemma \ref{l:I_q}(ii),
for every $\theta>0$,
$$
\left|\sum_{|c|>h^\delta}\sum_{k=1}^\infty 
k^{-4} S_{k}(c;\xi) I_{k,\ell}({}^tg c)\right|
\ll_{w,\ell,\theta} h^5 
\left(\sum_{k=1}^\infty 
k^{-5} |S_{k}(c;\xi)|\right)
\left(\sum_{|c|>h^\delta} |{}^tg c|^{-\theta} \right).
$$
By Lemma \ref{l:sq}(i), the first sum is finite.
For $\theta>4$, the second sum is estimated as 
$$
\sum_{|c|>h^\delta} |{}^tg c|^{-\theta} \ll \|g\|_E^\theta \sum_{|c|>h^\delta} |c|^{-\theta}
\ll_\theta \|g\|_E^\theta h^{-(\theta-4)\delta}.
$$
Hence,
\begin{equation}
\label{eq:step1}
\left|\sum_{|c|>h^\delta}\sum_{k=1}^\infty 
k^{-4} S_{k}(c;\xi) I_{k,\ell}({}^tg c)\right|
\ll_{w,\ell,\theta} \|g\|_E^\theta h^{5-(\theta-4)\delta}.
\end{equation}

Next, we estimate the terms with $0<|c|\le h^\delta$.
By Lemma \ref{l:I_q}(i)(iii), for every $\epsilon\in (0,1/2)$,
$$
\left|\sum_{0<|c|\le h^\delta}\sum_{k=1}^\infty 
k^{-4} S_{k}(c;\xi) I_{k,\ell}({}^tg c)\right|
\ll_{w,\ell,\epsilon} h^{3+2\epsilon} 
\left(\sum_{k=1}^{h} k^{-3 }|S_{k}(c;\xi)| \right)
\left(\sum_{0<|c|\le h^\delta} |{}^tgc|^{-(1-\epsilon)}\right)
.
$$
It follows from Lemma \ref{l:sq}(ii),
$$
\sum_{k=1}^{h} k^{-3 }|S_{k}(c;\xi)|\ll_{\ell,\epsilon}  h^{1/2+\epsilon},
$$
and 
$$
\sum_{0<|c|\le h^\delta} |{}^tgc|^{-(1-\epsilon)}\ll \|g\|_E\cdot \sum_{0<|c|\le h^\delta} |c|^{-(1-\epsilon)} \ll \|g\|_E h^{(3+\epsilon)\delta}.
$$
Hence, we conclude that for every $\epsilon>0$,
\begin{equation}
\label{eq:step2}
\left|\sum_{0<|c|\le h^\delta}\sum_{k=1}^\infty 
k^{-4} S_{k}(c;\xi) I_{k,\ell}({}^tg c)\right|
\ll_{w,\ell,\theta,\epsilon} \|g\|_E h^{7/2+3\delta+\epsilon}.
\end{equation}
Combining \eqref{eq:step1} and \eqref{eq:step2}, we deduce that
\begin{align*}
\hbox{N}_h(N,w_g; \xi)=&c_h h^{-2}\ell^{-4} \sum_{k=1}^\infty 
k^{-4} S_{k}(0;\xi) I_{k,\ell}(0)\\
&+
O_{w,\ell,\theta,\epsilon}\big(\|g\|_E^\theta h^{3-(\theta-4)\delta}+\|g\|_E h^{3/2+3\delta+\epsilon}\big).
\end{align*}
The last sum was already estimated as in the proof of Theorem \ref{th:heath-brown}.
\end{proof}

It will be convenient to interpret the local densities $\sigma_\infty$ and $\sigma_f$ group-theoretically, namely, in terms of the Tamagawa measures for the group $\sf G$.
We refer to \cite[Ch.~2]{w} or  \cite[Ch.~5]{vosk} 
for basic properties of the Tamagawa measures.
Let us fix a nowhere zero regular rational differential for of top degree on $\sf G$
(this form is known to be unique up to a constant factor).
Integration with respect to this form defines Haar measures  
$\tau_\infty$ and $\tau_q$
on $G_\infty:={\sf G}(\RR)$ and $G_q:={\sf  G}(\QQ_q)$ respectively. 
While the local Tamagawa measures are only unique up to constant factors,
remarkably their product is canonical.
In fact, according to the Tamagawa Volume Formula, since $\sf G$ is simply connected
\begin{equation}\label{eq:tamagawa}
\tau_\infty(G_\infty/\Gamma)\cdot {\prod}_{q\hbox{\tiny -prime}} \tau_q({\sf G}(\ZZ_q))=1.
\end{equation}

It will be convenient to define the Tamagawa measures using that $\sf G$ is a fiber of the norm map $N$. Let $\omega_0$ be the standard one-form on $\mathbb{A}^1$,
and $\omega$ be the standard form of top degree on ${\sf D}$.
Then there exists a form $\eta$ of degree three on $\sf D$ such that $\eta\wedge N^*(\omega_0)=\omega$. We denote by $\omega_y$ the restriction of $\eta$ to the fiber
$N^{-1}(y)$ for $y\ne 0$, and by $\tau_q^{(y)}$ the corresponding measures 
supports on the fibers $N^{-1}(y)(\Q_q)$. In particular, $\tau_q:=\tau_q^{(1)}$
defines a Tamagawa measures on ${\sf G}(\Q_q)$.
We note that the fiber measures 
are uniquely defined by the disintegration formula
\begin{equation}\label{eq:disintegration}
\int_{{\sf D}(\QQ_q)}\phi(N(x))\psi(x)dx=\int_{\QQ_q\backslash \{0\}} \phi(y) \left(\int_{N^{-1}(y)} \psi\, d\tau_q^{(y)}\right) dy,
\end{equation}
where $\phi$ and $\psi$ are compactly supported locally constant function.

Previously, we used the measures $m_\infty$ and $m_q$ 
that are normalized as 
$$
m_\infty(G_\infty/\Gamma)=1\quad\hbox{and}\quad m_q({\sf G}(\ZZ_q))=1.
$$
In view of \eqref{eq:tamagawa}, these measures can be expressed in terms of 
the Tamagawa measures as 
\begin{equation}\label{eq:mtau}
m_\infty=\tau_\infty(G_\infty/\Gamma)^{-1}\tau_\infty\quad\hbox{and}\quad
m_q=\tau_q({\sf G}(\ZZ_q))^{-1}\tau_q.
\end{equation}

\begin{lem}\label{l:local}
\begin{enumerate}
\item[(i)] Let $p$ be a prime, $\ell\in \N$ coprime to $p$, and $\xi \in {\sf G}(\Z/\ell\Z)$. Then for $h=p^s$,
$$
\sigma_f(N,\xi,h)
= \tau_p({\sf G}(\ZZ_p))^{-1}h^{-2} \tau_p(B_h)\ell^4\tau_\infty(G_\infty/\Gamma_\ell)^{-1},
$$
where $B_h=\{g\in G_p:\, \|g\|_p\le h \}$.
\item[(ii)] For every $w\in C_c(G_\infty)$,
$$
\sigma_\infty(N,w)=\int_{G_\infty} w\, d\tau_\infty.
$$
\end{enumerate}
\end{lem}

\begin{proof}
The crucial connection between the local densities $\sigma_q$ and the Tamagawa measures $\tau_q$ is given by the following formula: 
\begin{equation} \label{eq:local0}
\tau_q({\sf G}(\ZZ_q)) =\lim_{e\to\infty} q^{-3e} |{\sf G}(\ZZ/q^{e}\ZZ)|.
\end{equation}

Let $q$ be a prime coprime to both $\ell$ and $p$. Then
\begin{align*}
\hbox{N}_h(N,q^e,q^{s_q},\xi)&=\big|\{x\,\hbox{mod}\, q^{e}:\,\, N(x)=h^2\,\hbox{mod}\, q^e\}\big|\\
&=\big|\{x\,\hbox{mod}\, q^{e}:\,\, N(x)=1\,\hbox{mod}\, q^e\}\big|\\
&=|{\sf G}(\ZZ/q^e\ZZ)|.
\end{align*}
Hence, it follows from \eqref{eq:local0} that in this case
\begin{equation}
\label{eq:tau_1}
\sigma_q(\xi,h)=\tau_q({\sf G}(\ZZ_q)).
\end{equation}

Let $q$ be a prime dividing $\ell=\prod_{r\hbox{\tiny -prime}} r^{s_r}$ (and coprime to $h=p^s$). 
We choose a residue $\bar h\,(\hbox{mod}\, p^e)$ such that $\bar h h=1 \,(\hbox{mod}\, q^e)$. Then for $e\ge s_q$,
\begin{align*}
\hbox{N}_h(N,q^e,q^{s_q},\xi)
&=
\big|\{x\,\hbox{mod}\, q^{e+s_q}:\,\, x=\xi \,\hbox{mod}\, q^{s_q},\, N(x)=h^2\,\hbox{mod}\, q^e\}\big|\\
&=
q^{4s_q}\big|\{y\,\hbox{mod}\, q^{e}:\,\, y=\xi \,\hbox{mod}\, q^{s_q},\, N(y)=h^2\,\hbox{mod}\, q^e\}\big|\\
&=q^{4s_q}
\big|\{y\,\hbox{mod}\, q^{e}:\,\, y=\bar h^2\xi \,\hbox{mod}\, q^{s_q},\, N(y)=1\,\hbox{mod}\, q^e\}\big|\\
&=q^{4s_q} \frac{|{\sf G}(\ZZ/q^e\ZZ)|}{|{\sf G}(\ZZ/q^{s_q}\ZZ)|}.
\end{align*}

Hence, by \eqref{eq:local0} as before,
\begin{equation}
\label{eq:tau_2}
\sigma_q(N,\xi,h)=q^{4s_q} |{\sf G}(\ZZ/q^{s_q}\ZZ)|^{-1} \tau_q({\sf G}(\ZZ_q)).
\end{equation}

Finally, we claim that
\begin{equation}
\label{eq:tau_3}
\sigma_p(N,\xi,h)=h^{-2} \tau_p(B_h).
\end{equation}
For $z\in \Q_p^\times$, let us consider the map $\Phi_z(x):=z^{-1}x$.
We observe that the map $\Phi$ transforms the fiber $N^{-1}(y)$ to $N^{-1}(z^{-2}y)$ in \eqref{eq:disintegration}, so that it follows from uniqueness of this integral decomposition that
\begin{equation}
\label{eq:change}
(\Phi_z)_*\big(\tau_p^{(y)}\big)=|z|_p^{2}\, \tau_p^{(z^{-2}y)}.
\end{equation}
For $r>0$ and $y\in \Q_p$, we set
$$
B_r(y):=\{x\in {\sf D}(\Q_p):\,\, N(x)=y,\, \|x\|_p\le r\}.
$$
Then
$$
\{x\in {\sf D}(\Z_p):\, N(x)=h^2 \,\hbox{mod}\, p^{e} \}={\bigsqcup}_{y\in h^2+p^e\Z_p} B_1(y),
$$
and it follows from \eqref{eq:disintegration} that
$$
\int_{h^2+p^e\ZZ_p} \tau_p^{(y)}(B_1(y))\,dy=p^{-4e}\big|\{x\in {\sf D}(\ZZ/p^e\ZZ):\, N(x)=h^2 \,\hbox{mod}\, p^{e} \}\big|.
$$
Hence,
$$
\sigma_p(N,\xi,h)=\lim_{e\to\infty} p^{e}\int_{h^2+p^e\ZZ_p} \tau_p^{(y)}(B_1(y))\,dy.
$$
If $e$ is sufficiently large, for every $y\in h^2+p^e\ZZ_p$ there exists $z$ such that 
$z^2=y$. Clearly, $|z|_p=|h|_p=h^{-1}$. Then by \eqref{eq:change},
$$
\tau_p^{(y)}(B_1(y))=\tau_p^{(y)}\Big(\Phi_z^{-1}\big(B_{|z|^{-2}_p}(1)\big)\Big)=|z|_p^{2}\,\tau_p^{(1)}\left(B_{|z|^{-2}_p}(1)\right)=h^{-2}\,\tau_p(B_h).
$$
This implies \eqref{eq:tau_3}.

Furthermore, $\Phi$  defines a bijection between
$$
B'_h:=\{x\in {\sf D}(\Q_p):\,\, N(x)=h^2,\, \|x\|_p\le 1\}=\{x\in {\sf D}(\ZZ_p):\, N(x)=h^2\}
$$
and
$$
B_h=\{x\in {\sf D}(\QQ_p):\,\, N(x)=1,\,\|x\|_p\le h\}.
$$
Hence, it follows that
$$
\tau_p(B_h)=\tau_p^{(1)}\big(\Phi(B_h')\big)=|h|_p^{-2}\,\tau_p^{(h^2)}(B'_h)=h^{2}\,\tau_p^{(h^2)}(B'_h).
$$
Furthermore, it follows from \eqref{eq:disintegration} that
$$
\int_{h^2+p^e\ZZ_p} \tau_p^{(y)}(B'_h)\,dy=p^{-4e}\big|\{x\in (\ZZ/p^e\ZZ)^4:\, N(x)=h^2 \,\hbox{mod}\, p^{e} \}\big|.
$$
Hence,
$$
\sigma_p(N,\xi,h)=\lim_{e\to\infty} p^{e}\int_{h^2+p^e\ZZ_p} \tau_p^{y}(B_h)dy 
=\tau_p^{h^2}(B'_h)=h^{-2} \tau_p(B_h).
$$
This proves \eqref{eq:tau_3}.

Combining \eqref{eq:tau_1}, \eqref{eq:tau_2} and \eqref{eq:tau_3}, we deduce that
\begin{align*}
\sigma_f(N,\xi,h)= \ell^4  |{\sf G}(\ZZ/\ell\ZZ)|^{-1} \left({\prod}_{q\ne p} \tau_q({\sf G}(\ZZ_q))\right) h^{-2}
\tau_p(B_h)
\end{align*}
Furthermore, by the Tamagawa Formula \eqref{eq:tamagawa}
\begin{align*}
\sigma_f(N,\xi,h) &=\ell^4 |G(\ZZ/\ell\ZZ)|^{-1} \tau_p({\sf G}(\ZZ_p))^{-1 } \tau_\infty(G_\infty/\Gamma)^{-1} h^{-2}\tau_p(B_h)\\
&= \tau_p({\sf G}(\ZZ_p))^{-1}h^{-2} \tau_p(B_h)\ell^4\tau_\infty(G_\infty/\Gamma_\ell)^{-1}.
\end{align*}
This proves the first formula.

It follows from the disintegration formula \eqref{eq:disintegration} that 
\begin{align*}
\sigma_\infty(N,w)&=\lim_{\epsilon\to 0^+} (2\epsilon)^{-1}\int_{|N(x)-1|\le \epsilon} w(x)\, dx\\
&=\lim_{\epsilon\to 0^+} (2\epsilon)^{-1}\int_{1-\epsilon}^{1+\epsilon} \left(\int_{N^{-1}(y)} w\, d\tau^{(y)}_\infty\right)\, dy=
\int_{N^{-1}(1)} w\, d\tau^{(1)}_\infty,
\end{align*}
which proves the second equality.
\end{proof}

From now on we fix prime $p$ and $\ell\in\N$ coprime to $p$. 
We always choose $h$ to be of the form $h=p^s$.

Our next goal is to show 
that the counting function that we studied can be interpreted in terms of 
the averaging operators $\rho_{p,\ell}(\beta_h)$. Let $w\in C_c(D_\infty)$, and  
for $x\in G_\infty$, we also set 
$$
w_x(g):=w(x^{-1}g).
$$
Let $\chi$ denote the the characteristic function of the subset ${\sf G}(\ZZ_p)$.
We introduce a function $\phi_w$ defined by 
\begin{equation}\label{eq:phii}
\phi_w(g_\infty,g_p):=\sum_{\gamma\in \Gamma_{p,\ell}} w(g_\infty\gamma)\chi(g_p\gamma) \quad \hbox{for $(g_\infty,g_p)\in G_\infty\times G_p.$}
\end{equation}
This defines the function on $X_{p,\ell}=(G_\infty\times G_p)/\Gamma_{p,\ell}$.
The invariant probability measure $\mu_{p,\ell}$ on $X_{p,\ell}$
is defined as 
$$
\int_{X_{p,\ell}}\left(\sum_{\gamma\in \Gamma_{p,\ell}} f(g\gamma)\right)\, d\mu_{p,\ell}(g\Gamma_{p,\ell}) =|\Gamma:\Gamma_{\ell}|^{-1}\int_{G_\infty\times G_p}
f\, d(m_\infty\times m_p)
$$
for $f\in C_c(G_\infty\times G_p)$. Indeed, if $F$ is a fundamental domain for $\Gamma_\ell $ in $G_\infty$, then $F\times {\sf G}(\Z_p)$
is a fundamental domain for $\Gamma_{p,\ell}$ in $G_\infty\times G_p$.
Since $m_\infty(G_\infty/\Gamma)=1$ and $m_p({\sf G}(\Z_p))=1$, 
the above formula indeed defines the invariant probability measure on 
$X_{p,\ell}$. In particular, it follows that
\begin{align*}
\int_{X_{p,\ell}} \phi_w \, d\mu_{p,\ell}&=
|\Gamma:\Gamma_{\ell}|^{-1}\left(\int_{G_\infty} w \, dm_\infty\right) \left( \int_{G_p}
\chi \, dm_p\right)\\ 
&=|{\sf G}(\ZZ/\ell\ZZ)|^{-1}\int_{G_\infty} w \, dm_\infty .
\end{align*}
Taking \eqref{eq:mtau} and Lemma \ref{l:local} into account, we obtain:
\begin{align}
\sigma_\infty(N,w)\sigma_f(N,\xi,h)
&= \left(\int_{G_\infty} w\, d\tau_\infty\right)\tau_p({\sf G}(\ZZ_p))^{-1}h^{-2} \tau_p(B_h)\ell^4\tau_\infty(G_\infty/\Gamma_\ell)^{-1} \nonumber\\
&= h^{-2} \ell^4 \left(\int_{X_{p,\ell}} \phi_w \, d\mu_{p,\ell}\right)m_p(B_h).\label{eq:i1}
\end{align}

For $u\in {\sf G}(\ZZ_p)$, we obtain:
\begin{align*}
\int_{B_h} \phi_w\big(b^{-1}(x^{-1},u)\big)\, dm_p(b)&=\int_{B_h} \left(\sum_{\gamma\in \Gamma_{p,\ell}}w(x^{-1}\gamma)\chi(b^{-1}u\gamma)\right)\,dm_p(g)\\
&= \sum_{\gamma\in \Gamma_{p,\ell}} w(x^{-1}\gamma) m_p(u\gamma {\sf G}(\ZZ_p)\cap B_h)\\
&= \sum_{\gamma\in \Gamma_{p,\ell}\cap B_h} w_x(\gamma).
\end{align*}
Since
\begin{align*}
\Gamma_{p,\ell}\cap B_h &=\big\{x\in \Lambda[1/p]:\,\, N(x)=1,\, x= I\,(\hbox{mod}\, \ell),\, \|x\|_p\le h  \big\}\\
&=\big\{h^{-1}y:\, y\in\Lambda,\, N(y)=h^2,\, y= hI\,(\hbox{mod}\, \ell)\big\},
\end{align*}
we conclude that for any $u\in {\sf G}(\ZZ_p)$,
\begin{equation}\label{eq:i2}
\int_{B_h} \phi_w\big(b^{-1}(x^{-1},u)\big)\, dm_p(b)
=\hbox{N}_h\left(N,w_x,h I \,\hbox{mod}\, \ell\right).
\end{equation}

We will also the following result about integrability of the function $\phi_w$:

\begin{lem}\label{l:lp}
There exists $\epsilon_0>0$ such that if $\hbox{\rm supp}(w)\subset B(e,\epsilon_0)$,
then $\phi_w\in L^r(X_{p,\ell})$ for all $r\in [1,\infty)$.
\end{lem}

\begin{proof}
Without loss of generality, we may assume that $w\ge 0$.

It is sufficient to show that $\phi_w\in L^r(X_{p,\ell})$ fo every $r\in \N$. We obtain: 
\begin{align*}
\|\phi_w\|_{L^r(X_{p,\ell})}^r&=\int_{X_{p,\ell}} \phi_w(x)^r \, d\mu_{p,\ell}(x)\\
&=\int_{(G_\infty\times G_S)/\Gamma_{p,\ell}}\left( \sum_{\gamma_1,\ldots,\gamma_r\in\Gamma_{p,\ell}} \phi_w(g\gamma_1)\cdots \phi_w(g\gamma_r) \right) d\mu_{p,\ell}(g\Gamma_{p,\ell}) \\
&=\int_{(G_\infty\times G_S)/\Gamma_{p,\ell}}\left( \sum_{\gamma_1,\ldots,\gamma_r\in\Gamma_{p,\ell}} \phi_w(g\gamma_1)\phi_w(g\gamma_1\gamma_2)\cdots \phi_w(g\gamma_1\gamma_r) \right) d\mu_{p,\ell}(g\Gamma_{p,\ell}) \\
&=\int_{G_\infty\times G_S}\left( \sum_{\gamma_2,\ldots,\gamma_r\in\Gamma_{p,\ell}} \phi_w(g)\phi_w(g\gamma_2)\cdots \phi_w(g\gamma_r) \right) d(m_\infty\times m_p)(g).
\end{align*}
We observe that the product is zero unless
$$
\gamma_i\in \hbox{supp}(\phi_w)^{-1}\hbox{supp}(\phi_w)\subset B(e,\epsilon_0)^{-1} B(e,\epsilon_0) \times {\sf G}(\Z_p).
$$
When $\epsilon_0$ is chosen to be sufficiently small, this implies that $\gamma_i=e$.
Hence, we conclude that 
\begin{align*}
\|\phi_w\|_r^r=\int_{G_\infty\times G_S}\phi_w^r\, d(m_\infty\times m_p)=\|w\|^r_{L^r(G_\infty)}<\infty.
\end{align*}
\end{proof}

The following proposition verifies Theorem \ref{th:sl2} for the class of functions $\phi_w$. While this class of functions is quite "sparse", we will eventually show 
that this can be used to derive this estimate for general functions.

\begin{prop}\label{p:norm}
Let $w\in C_c^\infty(D_\infty)$ be as in Lemma \ref{l:lp}.
Then for every $\epsilon>0$,	
$$
\left\|\rho_{p,\ell}(\beta_h)\phi_w-\int_{X_{p,\ell}}\phi_w\, d\mu_{p,\ell}\right\|_{L^2(X_{p,\ell})}\ll_{w,p,\ell,\epsilon} m_p(B_h)^{-\sigma+\epsilon},
$$
where $\sigma=1/4$ if ${\sf G}$ is a anisotropic over $\Q$ and $\sigma=1/16$ otherwise.
\end{prop}

\begin{proof}
Taking \eqref{eq:i1} and \eqref{eq:i2} into account,
Theorem \ref{th:heath-brown2} 
can be restated as follows:
for every $g\in G_\infty$, $u\in {\sf G}(\Z_p)$, $\delta>0$, $\theta>4$, and $\epsilon>0$,
\begin{align}\label{eq:mmm}
\rho_{p,\ell}(\beta_h) \phi_w(g^{-1},u)=&\int_{X_{p,\ell}}\phi_w\, d\mu_{p,\ell}\\
&\;+O_{p,w,\ell,\delta,\theta,\epsilon}\Big(\|g\|_E^\theta h^{1-(\theta-4)\delta}+\|g\|_E h^{-1/2+3\delta+\epsilon} \Big).\nonumber
\end{align}
We shall use this estimate to prove the proposition.

First, we note that 
that \eqref{eq:mmm} implies 
for every compact $Q\subset G_\infty$,
\begin{align*}
\left\|\rho_{p,\ell}(\beta_h) \phi_w((\cdot)^{-1},u)-\int_{X_{p,\ell}}\phi_w\, d\mu_{p,\ell}\right\|_{L^2(Q)}
\ll_{p,w,\ell,Q,\epsilon} h^{-1/2+\epsilon}
\end{align*}
for all $\epsilon>0$. If $\sf G$ is anisotropic over $\QQ$, 
then the lattice $\Gamma_{\ell}$ is cocompact in $G_\infty$.
Hence, we can choose a compact $Q$ such that $Q^{-1}$ surjects  onto $G_\infty/\Gamma_\ell$, and the second part of the corollary follows.
Then $Q^{-1}\times {\sf G}(\Z_p)$ surjects onto $(G_\infty\times G_p)/\Gamma_{p,\ell}$.
Hence, it follows from the above estimate that 
\begin{align*}
\left\|\rho_{p,\ell}(\beta_h) \phi_w-\int_{X_{p,\ell}}\phi_w\, d\mu_{p,\ell}\right\|_{L^2(X_{p,\ell})}
\ll_{p,w,\ell,\epsilon} h^{-1/2+\epsilon}\ll_p m_p(B_h)^{-1/4+\epsilon}.
\end{align*}
This proves the proposition when $\sf G$ is anisotropic over $\QQ$.

Now we consider the case when $\sf G$ is isotropic over $\QQ$.
We observe that if $Q$ is a subset of $G_\infty$ such that $Q^{-1}$ surjects  onto $G_\infty/\Gamma_\ell$, then as before
\begin{align*}
\left\|\rho_{p,\ell}(\beta_h) \phi_w-\int_{X_{p,\ell}}\phi_w\, d\mu_{p,\ell}\right\|_{L^2(X_{p,\ell})}
\le \left\|\rho_{p,\ell}(\beta_h) \phi_w((\cdot)^{-1},e)-\int_{X_{p,\ell}}\phi_w\, d\mu_{p,\ell}\right\|_{L^2(Q)}.
\end{align*}
It follows from the theory of Siegel sets that such $Q$ can be chosen of the form
$$
Q=Q_0\cup Q_1\cup\cdots \cup Q_s,
$$ 
where $Q_0$ is compact, and for $i\ge 1$,
\begin{equation}
\label{eq:q_i}
Q_i^{-1}:=\{ka(t)ng_i:\, k\in K,t\ge 0, n\in N_0\}.
\end{equation}
Here $K=\hbox{SO}(2)$, $a(t)=\hbox{diag}(e^t,e^{-t})$,
$N_0$ is a compact subset of the upper triangular unipotent group, and $g_i\in {\sf G}(\Q)$. 

Let us consider the case when $Q$ is given by \eqref{eq:q_i}.
We note that the case of the union can be handled by using the triangle inequality.
We set 
$$
Q_{<R}:=\{t<\log R\}\quad\hbox{ and } \quad Q_{\ge R}:=\{t\ge \log R\}.
$$
We observe that for every $\theta>2$,
$$
\int_{Q_{<R}} \|g\|_E^\theta\, dm_\infty(g)\ll \int_{0}^{\log R} e^{(\theta-2)t}\,dt\ll R^{\theta-2},
$$
and
$$
\int_{Q_{<R}} \|g\|_E\, dm_\infty(g)\ll 1.
$$
Hence, it follows from \eqref{eq:mmm} that
$$
\left\|\rho_{p,\ell}(\beta_h) \phi_w((\cdot)^{-1},e)-\int_{X_{p,\ell}}\phi_w\, d\mu_{p,\ell}\right\|_{L^2(Q_{<R})}
\ll _{p,w,\ell,\delta,\theta,\epsilon} R^{\theta-2} h^{1-(\theta-4)\delta}+ h^{-1/2+3\delta+\epsilon}.
$$

To estimate the integral over $Q_{\ge R}$, we consider the set  
$$
\Omega_R:=\big(Q^{-1}_{\ge R}\times {\sf G}(\Z_p)\big)\Gamma_{p,\ell}\subset X_{p,\ell},
$$
and $\omega_R$ denote the characteristic function of this set. Then
\begin{align*}
&\left\|\rho_{p,\ell}(\beta_h) \phi_w((\cdot)^{-1},e)-\int_{X_{p,\ell}}\phi_w\, d\mu_{p,\ell}\right\|_{L^2(Q_{\ge R})}\\
\ll &
\left\|\left(\rho_{p,\ell}(\beta_h) \phi_w -\int_{X_{p,\ell}}\phi_w\, d\mu_{p,\ell}\right) \omega_R\right\|_{L^2(X_{p,\ell})}.
\end{align*}
Using the H\"older inequality, we obtain that for every $r\ge 1$ and $s=(1-1/r)^{-1}$,
\begin{align*}
&\left\|\left(\rho_{p,\ell}(\beta_h) \phi_w -\int_{X_{p,\ell}}\phi_w\, d\mu_{p,\ell}\right) \omega_R\right\|_{L^2(X_{p,\ell})}^2\\
\le & \left\|\rho_{p,\ell}(\beta_h) \phi_w -\int_{X_{p,\ell}}\phi_w\, d\mu_{p,\ell}\right\|_{L^{2r}(X_{p,\ell})}^2\cdot  \left\|\omega_R\right\|_{L^{2s}(X_{p,\ell})}^2.
\end{align*}
It follows from Jensen inequality that 
the operator $\rho_p(\beta_h): L^{2r}(X_{p,\ell})\to L^{2r}(X_{p,\ell})$
is bounded and the corresponding norm satisfies $\|\rho_p(\beta_h)\|\le 1$.
Hence, since $\phi_w\in  L^{2r}(X_{p,\ell})$ by Lemma \ref{l:lp}, we conclude that 
\begin{align*}
\left\|\rho_{p,\ell}(\beta_h) \phi_w((\cdot)^{-1},e)-\int_{X_{p,\ell}}\phi_w\, d\mu_{p,\ell}\right\|_{L^2(Q_{\ge R})}
&\ll_{w,p,\ell} \left\|\omega_R\right\|_{L^{2s}(X_{p,\ell})}\\
&=\mu_{p,\ell}(\Omega_R)^{1/(2s)}\\
&\le m_\infty(Q_{\ge R})^{1/(2s)}\ll_s R^{-1/s}.
\end{align*}
This implies that for every $\epsilon>0$, 
$$
\left\|\rho_{p,\ell}(\beta_h) \phi_w((\cdot)^{-1},e)-\int_{X_{p,\ell}}\phi_w\, d\mu_{p,\ell}\right\|_{L^2(Q_{\ge R})}
\ll_{w,p,\ell,\epsilon } R^{-1+\epsilon}.
$$
Ultimately, we conclude that 
\begin{align*}
\left\|\rho_{p,\ell}(\beta_h) \phi_w-\int_{X_{p,\ell}}\phi_w\, d\mu_{p,\ell}\right\|_{L^2(X_{p,\ell})}
\ll _{p,w,\ell,\delta,\theta,\epsilon} R^{\theta-2} h^{1-(\theta-4)\delta}+ h^{-1/2+3\delta+\epsilon}+R^{-1+\epsilon}.
\end{align*}
We choose  $R=h^{1/2-3\delta}$. Then for every $\epsilon>0$,
\begin{align*}
\left\|\rho_{p,\ell}(\beta_h) \phi_w-\int_{X_{p,\ell}}\phi_w\, d\mu_{p,\ell}\right\|_{L^2(X_{p,\ell})}
\ll_{w,p,\ell,\theta,\delta,\epsilon} h^{-\sigma+\epsilon},
\end{align*}
where $\sigma:=\min\big(-(1/2-3\delta)(\theta-2)+(\theta-4)\delta-1, 1/2-3\delta\big)$.
To optimise the error term, we choose $\delta=(\theta+1)/(8\theta-14)$.
Then as $\theta\to \infty$, we get $\sigma\to 1/8$. This implies the theorem.
\end{proof}

\begin{proof}[Proof of Theorem \ref{th:sl2}]
We recall that by Proposition \ref{p:norm}, for every $w\in C_c^\infty(D_\infty)$,
$$
\left\|\rho_{p,\ell}(\beta_h)\phi_w-\int_{X_{p,\ell}}\phi_w\, d\mu_{p,\ell}\right\|_{L^2(X_{p,\ell})}\ll_{w,p,\ell,\epsilon} m_p(B_h)^{-\sigma+\epsilon}
$$
for all $\epsilon>0$, where $\phi_w\in L^2(X_{p,\ell})$ is defined by \eqref{eq:phii}.
Let us also consider a family of functions $w_x(y):=w(x^{-1}y)$ with $x\in G_\infty$.
We observe that $\phi_{w_x}=\rho_{\infty,\ell}(x)(\phi_w)$ for the operator 
$\rho_{\infty,\ell}(x):L^2(X_{p,\ell})\to L^2(X_{p,\ell}): \phi\mapsto \phi\circ x^{-1}$.
Since $\rho_{\infty,\ell}(x)$ commutes with $\rho_{p,\ell}(\beta_h)$ and $\|\rho_{\infty,\ell}(x)\|=1$, we deduce 
that for every $x\in G_\infty$,
$$
\left\|\rho_{p,\ell}(\beta_h)\phi_{w_x}-\int_{X_{p,\ell}}\phi_{w_x}\, d\mu_{p,\ell}\right\|_{L^2(X_{p,\ell})}\ll_{w,p,\ell,\epsilon} m_p(B_h)^{-\sigma+\epsilon}.
$$
Let $\pi$ be an irreducible unitary representation of $G_p$ which is
discretely embedded in $\rho_{p,\ell}$.
Since the sets $B_h$ are ${\sf G}(\Z_p)$-invariant, $\pi(\beta_h)=0$.
Hence, we may assume that $\pi$ is spherical and 
denote by $F_\pi\in L^2(X_{p,\ell})$ the unique unit ${\sf G}(\Z_p)$-invariant
vector associated to $\pi$. Arguing exactly as in the proof of Theorem \ref{th:converse}, we deduce that 
$$
\|\pi(\beta_h)\| 
\ll_{w,p,\ell,\epsilon}  \left|\left<\phi_{w_x},F_\pi\right>\right|^{-1} m_p(B_h)^{-\sigma+\epsilon},
$$
provided that $\left<\phi_{w_x},F_\pi\right>\ne 0$. Moreover,
$$
\left<\phi_{w_x},F_\pi\right>=\int_{G_\infty} w(x^{-1}g)f(g)\, dm_\infty(g),
$$
where $f(g):= F_\pi(g\Gamma_{p,\ell})$.
Since $f$ is a non-zero function which is locally $L^2$-integrable,
it follows from the following general version of 
the Local Ergodic Theorem (Lemma \ref{l:local} below) that 
there exists $x\in G_\infty$ such that $\left<\phi_{w_x},F_\pi\right>\ne 0$.
Hence, we conclude that
\begin{equation}
\label{eq:s3}
\|\pi(\beta_h)\| 
\ll_{\pi,w,p,\ell,\epsilon}   m_p(B_h)^{-\sigma+\epsilon}
\end{equation}
for all $\epsilon>0$.
We refer to \cite[Ch.~2]{ggp} for the classification of the irreducible unitary
representations of $G_p\simeq \hbox{SL}_2(\Q_p)$. In particular, 
let us consider the complementary series representations $\pi_s\in \widehat G_p$ with $s\in (0,1)$. These representations a spherical, and 
we recall that the corresponding spherical function are estimated as 
\begin{equation}
\label{eq:s1}
\|g\|_p^{-(1-s)} \ll_{p,s}|\omega_{\pi_s}(g)|\ll_{p,s} \|g\|_p^{-(1-s)}\quad \hbox{for $g\in G_p$.}
\end{equation}
Using this bound, we deduce from \eqref{eq:norm_spher} that
$$
h^{-(1-s)} \ll_{p,s} \|\pi_s(\beta_h)\|\ll_{p,s} h^{-(1-s)}.
$$
Since also
$$
h^2\ll_p m_p(B_h)\ll_p h^2,
$$
we conclude that
\begin{equation}
\label{eq:s2}
m_p(B_h)^{-(1-s)/2} \ll_{p,s} \|\pi_s(\beta_h)\|\ll_{p,s} m_p(B_h)^{-(1-s)/2}.
\end{equation}
It will be important for us that the implicit constants in \eqref{eq:s1} and hence in \eqref{eq:s2} are uniformly bounded for $s\le s_0<1$.

Comparing \eqref{eq:s3} and \eqref{eq:s2} when $m_p(B_h)\to\infty$,
we deduce that if 
the complementary series representation $\pi_s$ is discretely embedded in 
$\rho_{p,\ell}$, then $s\le 1-2\sigma$.
Hence, we conclude that if $\pi_s\in \widehat G_p^{{\rm aut},0}$, then
$s\le 1/2$ if $\sf G$ is anisotropic over $\Q$, and 
$s\le 7/8$ if $\sf G$ is anisotropic over $\Q$.
The continuous component of the representations $\rho_{p,\ell}$ has been
described (in much greater generality) by Langlands \cite{L76}.
It follows from this description that the continuous component is tempered.
Moreover, it follows from the description of the unitary dual of 
$\hbox{SL}_2(\Q_p)$ (see, for instance, \cite[Ch.~2]{ggp}) that the only non-tempered irreducible unitary representations
are the complementary series $\pi_s$. Therefore,
$$
\rho_{p,\ell}= \left({\sum}_i \pi_{s_i}^{\oplus n_i}\right)  \oplus \rho_{p,\ell}',
$$
where $\rho_{p,\ell}'$ is a tempered representation and $s_i\le 1-2\sigma$ for all $i$.
Since for tempered representations the bound \eqref{eq:temp} holds, we conclude that 
$$
\|\rho_{p,\ell}(\beta_h)\|\ll_{p,\ell} m_p(B_h)^{-(1-s_{\rm\tiny max})/2},
$$
where $s_{\rm\tiny max}:=\max(s_i)$. This completes the proof of the theorem modulo Lemma \ref{l:local}.
\end{proof}

\begin{lem} \label{l:local}
There exists a collection of smooth non-negative compactly supported function $w^{(r)}$, $r\in (0,r_0)$, on $G_\infty$ such that for every locally $L^2$-integrable $f$ on $G_\infty$,
\begin{equation}
\label{eq:local}
\left(\int_{G_\infty} w^{(r)} \, dm_\infty\right)^{-1} \int_{G_\infty} w^{(r)} (x^{-1}g)f(g)\, dm_\infty(g)\longrightarrow f(g)\quad\hbox{as $r\to 0^+$,}
\end{equation}
 for almost all $x\in G_\infty$.
\end{lem}

\begin{proof} 
We fix a $G_\infty$-left-invariant Riemannian on ${\sf D}(\RR)$ and consider $$
w^{(r)}(g):=\phi\big(r^{-1} d(g,e)\big),$$ where $\phi$ is a smooth non-negative symmetric bump function at 0.
Moreover, we assume that $\phi$ is non-increasing on $\RR^+$.
We note that the claim of the lemma clearly holds for continuous functions.
	
We set 	
$$
A_r f(x):=\left(\int_{G_\infty} w^{(r)} \, dm_\infty\right)^{-1} \int_{G_\infty} w^{(r)} (x^{-1}g)f(g)\, dm_\infty(g),
$$
and define the corresponding maximal function 
$$
M_\phi f(x):=\sup_{r\in (0,r_0)} A_r|f|(x).
$$
If the maximal function satisfies the bound 
\begin{equation}
\label{eq:max0}
\|M_\phi f\|_2\ll \|f\|_2
\end{equation}
for every $L^2$-integrable $f$,
then by a standard argument one can extend \eqref{eq:local} from continuous functions to general
$L^2$-integrable functions. Hence, it remains to prove \eqref{eq:max0}.
In fact, it can be deduced from the classical maximal inequality for the operators
$$
Mf(x):= \sup_{r>0} m_\infty(B(x,r))^{-1} \int_{B(x,r)} |f(g)|\, dm_\infty(g),
$$
where $B(x,r)$ denotes the balls in $G_\infty$ with respect to the metric $d$.
We choose positive parameters $\alpha_i,r_i=O_\phi(1)$ so that 
$$
\phi\le {\sum}_i \alpha_i\, \chi_{B(0,r_i)}\quad\quad \hbox{and}\quad\quad  {\sum}_i\alpha_ir_i\ll_\phi 1.
$$
Then 
\begin{align*}
\int_{G_\infty} w^{(\epsilon)} (x^{-1}g)|f(g)|\, dm_\infty(g)&\le 
\sum_i \alpha_i \int_{B(x,r_i\epsilon)} |f(g)|dm_\infty(g)\\
&\le \left( \sum_i\alpha_i m_\infty(B(x,r_i\epsilon))\right) Mf(x)\\
&\ll_\phi Mf(x).
\end{align*}
Hence, \eqref{eq:max0} follows from the classical maximal inequality.
\end{proof}

\end{document}